\documentclass[11pt]{article}

\usepackage{fullpage}
\usepackage[final]{hyperref}
\usepackage[round]{natbib}
\usepackage{tabularx}
\usepackage{booktabs}
\usepackage{setspace}
\usepackage{multirow}
\usepackage{makecell}
\usepackage{graphicx}
\usepackage{delarray}
\usepackage{xcolor}
\usepackage{colortbl}
\usepackage{amssymb}
\usepackage{amsmath}
\usepackage{mathrsfs}
\usepackage{mathtools}
\usepackage{nicefrac}       
\usepackage{tikz}
\usepackage{xspace}
\xspaceaddexceptions{[]\{\}}
\usepackage{algorithm}
\usepackage{algorithmic}
\usepackage{mdframed}
\usepackage{bm}
\usepackage{tablefootnote}
\usepackage{threeparttable}
\allowdisplaybreaks

\newtheorem{theorem}{Theorem}
\newtheorem{lemma}{Lemma}
\newtheorem{corollary}{Corollary}

\newtheorem{definition}{Definition}
\newtheorem{assumption}{Assumption}

{\hspace*{\fill}$\Box$\par}
{\hspace*{\fill}$\Box$\par\vspace{4mm}}
\newenvironment{proofof}[1]{\smallskip\noindent{\bf Proof of #1.}}%
{\hspace*{\fill}$\Box$\par}

\hypersetup{
	colorlinks=true,       
	linkcolor=black,        
	citecolor=blue,        
	filecolor=magenta,     
	urlcolor=blue
}

\newcommand{\eat}[1]{}
\newcommand{\topic}[1]{\vspace{2mm}\noindent{{\bf #1:}}}

\definecolor{bgcolor}{rgb}{0.66,0.88,1.00}

\newcommand{\E}{{\mathbb{E}}}
\newcommand{\R}{{\mathbb{R}}}
\newcommand{\inner}[2]{\langle #1,#2 \rangle}

\newcommand{\innerB}[2]{\Big\langle #1,#2 \Big\rangle}
\newcommand{\innerb}[2]{\big\langle #1,#2 \big\rangle}

\newcommand{\ns}[1]{\| #1 \|^2}

\newcommand{\n}[1]{\| #1 \|}

\newcommand{\hx}{\widehat{x}}

\newcommand{\anita}{{\sf ANITA}\xspace}

\newcommand{\uxt}{\underline{x}_t}

\newcommand{\xt}{x_t}
\newcommand{\xtn}{x_{t+1}}
\newcommand{\wt}{w_t}
\newcommand{\wtn}{w_{t+1}}

\newcommand{\bxt}{\bar{x}_t}
\newcommand{\bxtn}{\bar{x}_{t+1}}

\newcommand{\tnabla}{\widetilde{\nabla}_t}
\newcommand{\etat}{\eta_{t}}

\begin{document}

\title{\bf ANITA: An Optimal Loopless Accelerated Variance-Reduced Gradient Method}

\author{ Zhize Li \\
				Carnegie Mellon University \\
	    	\texttt{zhizeli@cmu.edu}
}

\date{}

\maketitle

\begin{abstract}
In this paper, we propose a novel accelerated gradient method called \anita for solving the fundamental finite-sum optimization problems.
Concretely, we consider both general convex and strongly convex settings:
i) For general convex finite-sum problems, \anita improves previous state-of-the-art result given by Varag \citep{lan2019unified}.
In particular, for large-scale problems or the convergence error is not very small, i.e., $n \geq \frac{1}{\epsilon^2}$, \anita obtains the \emph{first} optimal result $O(n)$, matching the lower bound $\Omega(n)$ provided by \citet{woodworth2016tight}, while previous results are $O(n \log \frac{1}{\epsilon})$ of Varag \citep{lan2019unified} and $O(\frac{n}{\sqrt{\epsilon}})$ of Katyusha \citep{allen2017katyusha}.
ii) For strongly convex finite-sum problems, we also show that \anita can achieve the optimal convergence rate $O\big((n+\sqrt{\frac{nL}{\mu}})\log\frac{1}{\epsilon}\big)$ matching the lower bound $\Omega\big((n+\sqrt{\frac{nL}{\mu}})\log\frac{1}{\epsilon}\big)$ provided by \citet{lan2015optimal}.
Besides, \anita enjoys a simpler loopless algorithmic structure unlike previous accelerated algorithms such as Varag \citep{lan2019unified} and Katyusha \citep{allen2017katyusha} where they use double-loop structures. Moreover, we provide a novel \emph{dynamic multi-stage convergence analysis}, which is the key technical part for improving previous results to the optimal rates. 
We believe that our new theoretical rates and novel convergence analysis for the fundamental finite-sum problem will directly lead to key improvements for many other related problems, such as distributed/federated/decentralized optimization problems (e.g., \citealp{li2021canita}).
Finally, the numerical experiments show that \anita converges faster than the previous state-of-the-art Varag \citep{lan2019unified}, validating our theoretical results and confirming the practical superiority of \anita.
\end{abstract}

\section{Introduction}
\label{sec:intro}

In this paper, we consider the fundamental finite-sum problems of the form
\begin{equation}\label{eq:prob}
\min_{x\in \R^d}   f(x):= \frac{1}{n}\sum_{i=1}^n{f_i(x)},
\end{equation}
where $f:\R^d\to \R$ is a smooth and convex function.
We consider two settings in this paper, i) general convex setting ($\mu=0$); ii) strongly convex setting ($\mu>0$), where $\mu$ is the strongly convex parameter for $f(x)$, i.e.,
$f(x) - f(y) - \inner{\nabla f(y)}{x-y} \geq \frac{\mu}{2}\ns{x-y}$. Note that the case $\mu=0$ reduces to the standard convexity.
Also note that the strong convexity is only corresponding to the average function $f$, is not needed for these component functions $f_i$s.  

Finite-sum problem \eqref{eq:prob} captures the standard empirical risk minimization (ERM) problems in machine learning \citep{shai_book}.
There are $n$ data samples and $f_i$ denotes the loss associated with $i$-th data sample, and the goal is to minimize the loss over all data samples. 
This optimization problem has found a wide range of applications in machine learning, statistical inference, and image processing.
In recent years, there has been extensive research in designing gradient-type methods for solving this problem \eqref{eq:prob}. 
To measure the efficiency of algorithms for solving \eqref{eq:prob}, it is standard to bound the number of stochastic gradient computations for finding a suitable solution. 
In particular, our goal is to find a point $\hx \in \R^d$ such that $\E[f(\hx)-f(x^*)] \leq \epsilon$, where the expectation is with respect to the randomness inherent in the algorithm. 
We use the term  {\em $\epsilon$-approximate solution} to refer to such a point $\hx$, 
and use the term \emph{stochastic gradient complexity} to describe the convergence result (convergence rate) of algorithms.

Two of the most classical gradient-type algorithms are gradient descent (GD) and stochastic gradient descent (SGD) (e.g., \citealp{nemirovski1983problem, nesterov2014introductory, nemirovski2009robust, duchi2010composite, lan2012optimal, ghadimi2012optimal, hazan2019introduction}).  
However, GD requires to compute the full gradient over all $n$ data samples for each iteration ($x_{t+1} = x_t-\eta \frac{1}{n}\sum_{i=1}^{n}\nabla f_i(x_t)$) which is inefficient especially for large-scale machine learning problems where $n$ is very large.
Although SGD only needs to compute a single stochastic gradient (e.g., $\nabla f_i(x)$) for each iteration ($x_{t+1}=x_t-\eta \nabla f_i(x_t)$), it requires an additional bounded variance assumption for the stochastic gradients (i.e., $\exists \sigma>0$, $\E_{i}[\ns{\nabla f_i(x)-\nabla f(x)}]\leq \sigma^2$) since it does not compute the full gradients ($\nabla f(x)$, i.e., $\frac{1}{n}\sum_{i=1}^{n}\nabla f_i(x)$).
More importantly, for strongly convex problems, SGD only obtains a sublinear convergence rate $O(\frac{\sigma^2}{\mu\epsilon})$ rather than a linear rate $O(\cdot\log\frac{1}{\epsilon})$ achieved by GD.

To remedy the variance term $\E[\ns{\nabla f_i(x)-\nabla f(x)}]$ in SGD, the variance reduction technique has been proposed and it has been widely-used in many algorithms in recent years.
In particular, \citet{leroux2012stochastic,schmidt2017minimizing} propose the first variance-reduced algorithm called SAG and show that by incorporating new gradient estimators into SGD one can possibly achieve the linear convergence rate for strongly convex problems.
Then this variance reduction direction is followed by many works such as \citep{shalev2013stochastic, mairal2013optimization, johnson2013accelerating, defazio2014saga, mairal2015incremental, nguyen2017sarah}.
Particularly, SAG \citep{leroux2012stochastic} uses a biased gradient estimator while SAGA \citep{defazio2014saga} modifies it to an unbiased estimator and provides better convergence results.
\citet{johnson2013accelerating} propose a novel unbiased stochastic variance reduced gradient (SVRG) method which directly incorporates the full gradient term $\nabla f(x)$ into SGD.
More specifically, each epoch of SVRG starts with the computation of the full gradient
$\nabla f(\tilde x)$ at a snapshot point $\tilde x\in \R^n$ and then runs SGD for a fixed number of steps
using the modified stochastic gradient estimator
\begin{align} \label{eq:svrg}
\tnabla = \nabla f_{i}(x_{t}) - \nabla f_{i}(\tilde x) + \nabla f(\tilde x),
\end{align}
i.e., $x_{t+1}=x_t-\eta \tnabla$, where $i$ is randomly picked from $\{1,2,\dots,n\}$.
In particular, if each full gradient $\nabla f(\tilde x)$ (which requires $n$ stochastic gradient computations) at the snapshot point $\tilde x$ is reused for $n$ iterations (i.e., $\tilde{x}$ is changed after every $n$ iterations), then the amortized stochastic gradient computations for each iteration is the same as SGD.
Note that $\E[\tnabla]=\nabla f(x_t)$ is an unbiased estimator, and its variance $\E[\ns{\tnabla -\nabla f(x_t)}] \leq 4L\big(f(x_t)-f(x^*) + f(\tilde{x})-f(x^*)\big)$ is reduced as the algorithm converges $x_t,\tilde{x} \rightarrow x^*$, while the variance term is uncontrollable for plain SGD where $\tnabla = \nabla f_{i}(x_{t})$.
\citet{johnson2013accelerating} also show that SVRG obtains the linear convergence $O((n+\frac{L}{\mu})\log \frac{1}{\epsilon})$ which can be better than the sublinear convergence rate $O(\frac{\sigma^2}{\mu\epsilon})$ of plain SGD, for strongly convex problems.
The SVRG gradient estimator \eqref{eq:svrg} is adopted in many algorithms (e.g., \citealp{xiao2014proximal, allen2015improved,lei2016less, allen2016variance, reddi2016stochastic, reddi2016proximal, lei2017non, li2018simple, zhou2018stochastic, ge2019stabilized, kovalev2019don}) and also is used in our \anita. 

The aforementioned variance-reduced methods are not accelerated and hence they do not achieve the optimal convergence rates for convex finite-sum problem \eqref{eq:prob}.
See the non-accelerated variance-reduced algorithms listed in the first part of Table~\ref{tab:1}, i.e., SAG, SVRG, SAGA and SVRG\textsuperscript{++}, they do not achieve the accelerated rates, i.e., $\frac{L}{\mu}$ vs.\ $\sqrt{\frac{L}{\mu}}$ (strongly convex case) and $\frac{L}{\epsilon}$ vs.\ $\sqrt{\frac{L}{\epsilon}}$ (general convex case). Note that we do not list the SCSG \citep{lei2016less} and SARAH \citep{nguyen2017sarah} in Table~\ref{tab:1} since SCSG requires an additional bounded variance assumption (without this assumption, its result is the same as SVRG and SAGA) and SARAH uses $\E[\ns{\nabla f(\hx)}]\leq \epsilon$ as the convergence criterion which can not be directly converted to $\E[f(\hx)-f(x^*)]\leq \epsilon$. SARAH is usually used for solving nonconvex problems where the convergence criterion is typically the norm of gradient (e.g., \citealp{fang2018spider,wang2018spiderboost,pham2019proxsarah,li2019ssrgd,li2021page}).
Also both SCSG and SARAH are non-accelerated methods and thus do not achieve the optimal convergence results. 
Therefore, much recent research effort has been devoted to the design of accelerated gradient methods (e.g., \citealp{nesterov2014introductory, beck2009fast, lan2012optimal, allen2014linear, su2014differential, lin2015universal, allen2017katyusha, lan2018random, lan2019unified, li2020fast,li2020acceleration}). 
As can be seen from Table~\ref{tab:1}, for strongly convex finite-sum problems, existing accelerated methods such as RPDG \citep{lan2015optimal}, Katyusha \citep{allen2017katyusha}, Varag \citep{lan2019unified} and our \anita are optimal since their convergence results are $O\big(\big(n+\sqrt{\frac{nL}{\mu}}\big)\log \frac{1}{\epsilon}\big)$ matching the lower bound $\Omega\big(\big(n+\sqrt{\frac{nL}{\mu}}\big)\log \frac{1}{\epsilon}\big)$ given by \citet{lan2015optimal}.

However, for general (non-strongly) convex finite-sum problems, all previous accelerated methods do not achieve the optimal convergence result. 
In particular, Varag \citep{lan2019unified} obtains the current best result $O\big(n\min\{\log\frac{1}{\epsilon}, \log n\}+\sqrt{\frac{nL}{\epsilon}}\big)$, while the lower bound in this general convex case is $\Omega\big(n+\sqrt{\frac{nL}{\epsilon}}\big)$ provided by \citet{woodworth2016tight}. 
More importantly, for large-scale problems where the number of data samples $n$ is very large, or the convergence error $\epsilon$ is not very small, then the convergence result of Varag is $O(n\log \frac{1}{\epsilon})$ which is not optimal since the lower bound is $\Omega(n)$ (see Table~\ref{tab:2}).
Note that the case of large-scale problems or the case of moderate convergence error often exists in machine learning applications.
We show that our \anita takes an important step towards the ultimate limit of accelerated methods and it is the first algorithm to achieve the optimal convergence rate $O(n)$ in this case matching the lower bound $\Omega(n)$.
See Tables~\ref{tab:1} and \ref{tab:2} for more details.

\begin{table}[t]
		\centering
		\caption{Convergence rates for finding an $\epsilon$-approximate solution $\E[f(\hx)-f(x^*)]\leq \epsilon$ of \eqref{eq:prob}}
		\label{tab:1}
		\vspace{1mm}
		\scriptsize
		\renewcommand{\arraystretch}{1.4}
		\begin{threeparttable}
			\begin{spacing}{1.4}
				\begin{tabular}{|c|c|c|c|}
					\hline
					\bf Algorithms
					& \bf $\mu$-strongly convex
					& \bf General convex
					& \bf \makecell{Loopless \\ \scriptsize{(Simple)}} \\
					\hline

					GD 
					& $O\left(\frac{nL}{\mu}\log\frac{1}{\epsilon}\right)$ 
					& $O\left(\frac{nL}{\epsilon}\right)$
					& Yes \\
					\hline
					
					\makecell{Nesterov's accelerated GD\\ 
						\citep{nesterov83, nesterov2014introductory}}
					& $O\left(n\sqrt{\frac{L}{\mu}}\log\frac{1}{\epsilon}\right)$ 
					& $O\left(n\sqrt{\frac{L}{\epsilon}}\right)$
					& Yes \\
					\hline
					
					\makecell{SAG~\citep{leroux2012stochastic}}
					& \hspace{-1.5mm}$O\Big(\big(n+n^2\lfloor\frac{L}{n \mu}\rfloor\big)\log\frac{1}{\epsilon}\Big)$ \hspace{-2.5mm}
					& ---
					& Yes \\
					\hline
					
					\makecell{SVRG~\citep{johnson2013accelerating}}
					& $O\left(\big(n+\frac{L}{\mu}\big)\log\frac{1}{\epsilon}\right)$ 
					& ---
					& No \\
					\hline
					
					\makecell{SAGA~\citep{defazio2014saga}}
					& $O\left(\big(n+\frac{L}{\mu}\big)\log\frac{1}{\epsilon}\right)$ 
					& $O\left(\frac{n+L}{\epsilon}\right)$
					& Yes \\
					\hline
					
					SVRG\textsuperscript{++}~\citep{allen2015improved}
					& ---
					& $O\left(n\log\frac{1}{\epsilon}+\frac{L}{\epsilon}\right)$
					& No \\
					\hline
					
					\makecell{RPDG~\citep{lan2015optimal}}
					& $O\left(\big(n+\sqrt{\frac{nL}{\mu}}\big)\log\frac{1}{\epsilon}\right)$ 
					& $O\left(\big(n+\sqrt{\frac{nL}{\epsilon}}\big)\log\frac{1}{\epsilon}\right)$ \tnote{1}
					& Yes \\
					\hline
					
					\makecell{Catalyst~\citep{lin2015universal}}
					&\hspace{-2mm}$O\left(\big(n+\sqrt{\frac{nL}{\mu}}\big)\log\frac{1}{\epsilon}\right)$ \tnote{1}\hspace{-2mm}
					& $O\left(\big(n+\sqrt{\frac{nL}{\epsilon}}\big)\log^2\frac{1}{\epsilon}\right)$ \tnote{1}
					& No \\
					\hline
					
					\makecell{Katyusha~\citep{allen2017katyusha}}
					& $O\left(\big(n+\sqrt{\frac{nL}{\mu}}\big)\log\frac{1}{\epsilon}\right)$ 
					& $O\left(n\log\frac{1}{\epsilon}+\sqrt{\frac{nL}{\epsilon}}\right)$ \tnote{1}
					& No \\
					\hline
					
					\makecell{Katyusha\textsuperscript{ns}~\citep{allen2017katyusha}}
					& ---
					& $O\left(\frac{n}{\sqrt{\epsilon}}+\sqrt{\frac{nL}{\epsilon}}\right)$ 
					& No \\
					\hline
					
					\makecell{Varag~\citep{lan2019unified}}
					& $O\left(\big(n+\sqrt{\frac{nL}{\mu}}\big)\log\frac{1}{\epsilon}\right)$ 
					& $O\Big(n\min\big\{\log\frac{1}{\epsilon},\ \log n\big\}+\sqrt{\frac{nL}{\epsilon}}\Big)$ 
					& No \\
					\hline

					 \gape{\multirow{2}{*}{\makecell{\anita (this paper)}}}
					&\cellcolor{bgcolor}  $O\left(\big(n+\sqrt{\frac{nL}{\mu}}\big)\log\frac{1}{\epsilon}\right)$
					&\hspace{-2.5mm} \cellcolor{bgcolor} $O\Big(n\min\big\{1+\log\frac{1}{\epsilon \sqrt{n}},\ \log\sqrt{n}\big\}+\sqrt{\frac{nL}{\epsilon}}\Big)$ \hspace{-3mm}
					& \cellcolor{bgcolor} Yes \\
					\cline{2-4}	
					
					&\cellcolor{bgcolor}
					$O\left(\big(n+\sqrt{\frac{nL}{\mu}}\big)\log\frac{1}{\epsilon}\right)$
					& \cellcolor{bgcolor} $O\left(n+\sqrt{\frac{nL}{\epsilon}}\right)$ 
					\tnote{2}
					& \cellcolor{bgcolor} Yes \\
					\hline	
					
					\makecell{Lower bound}
					& \makecell{$\Omega\left(\big(n+\sqrt{\frac{nL}{\mu}}\big)\log\frac{1}{\epsilon}\right)$\\ \citep{lan2015optimal}}
					& \makecell{$\Omega\left(n+\sqrt{\frac{nL}{\epsilon}}\right)$\\
						\citep{woodworth2016tight}} 
					& --- \\
					\hline
					
				\end{tabular}
			\end{spacing}
			\vspace{0mm}
			\begin{tablenotes}
				\scriptsize
				\item[1] These gradient complexity bounds are obtained via indirect approaches, i.e., by adding strongly convex perturbation.
				\item[2] \anita can achieve this optimal result for a very wide range of $\epsilon$, i.e., $\epsilon \in (0,\frac{L}{n\log^2\sqrt{n}}]\cup [\frac{1}{\sqrt{n}},+\infty)$
				or the number of data samples $n\in (0, \frac{L}{\epsilon\log^2\sqrt{n}}] \cup [\frac{1}{\epsilon^2}, +\infty)$	
				(see Table \ref{tab:2} for more details). Note that the term $\min\{\log\frac{1}{\epsilon},\ \log n\}$ in Varag~\citep{lan2019unified} cannot be removed regardless of the value of $\epsilon$ or $n$. 
				Thus \anita is the first accelerated algorithm that can exactly achieve the optimal convergence result.
			\end{tablenotes}
	\end{threeparttable}	
\end{table}
\begin{table}[!h]
	\centering
	\caption{Direct accelerated stochastic algorithms for  \emph{general convex setting} wrt. $\epsilon$} 
	\label{tab:2}
	\vspace{1mm}
	\scriptsize
	\renewcommand{\arraystretch}{1.4}
	\begin{spacing}{1.4}
		\begin{tabular}{|c|c|c|c|c|}
			\hline
			\multirow{2}{*}{Algorithms}
			&\multicolumn{4}{c|}{
				\makecell{
					The convergence error ($\E[f(\hx)-f(x^*)]\leq \epsilon$):  large $\epsilon$ $\longrightarrow$ small $\epsilon$ \\
					(or the number of data samples: large $n$ $\longrightarrow$ small $n$) }}\\
			\cline{2-5}

			&  \makecell{$\epsilon\geq \frac{1}{\sqrt{n}}$\\
				(or $n\geq\frac{1}{\epsilon^2}$)}
			
			&  \makecell{$\frac{1}{\sqrt{n}} > \epsilon\geq \frac{1}{n}$\\
				(or $\frac{1}{\epsilon^2}>n\geq \frac{1}{\epsilon}$)}
			&  \makecell{$\frac{1}{n} > \epsilon\geq \frac{L}{n\log^2\sqrt{n}}$\\
				(or $\frac{1}{\epsilon}>n\geq\frac{L}{\epsilon\log^2\sqrt{n}}$)}
			&  \makecell{$ \frac{L}{n\log^2\sqrt{n}} > \epsilon $\\
				(or $\frac{L}{\epsilon\log^2\sqrt{n}}>n$)} \\
			\hline
			
			\makecell{Katyusha\textsuperscript{ns}~\citep{allen2017katyusha}}
			& $O\left(\frac{n}{\sqrt{\epsilon}}\right)$
			& $O\left(\frac{n}{\sqrt{\epsilon}}\right)$
			& $O\left(\frac{n}{\sqrt{\epsilon}}\right)$
			& $O\Big(\frac{n}{\sqrt{\epsilon}}+\sqrt{\frac{nL}{\epsilon}}\Big)$  \\
			\hline
			
			\makecell{Varag~\citep{lan2019unified}}
			& $O\left(n\log\frac{1}{\epsilon}\right)$ 
			& $O\left(n\log\frac{1}{\epsilon}\right)$ 
			& $O\left(n\log n\right)$ 
			&  $O\Big(\sqrt{\frac{nL}{\epsilon}}\Big)$ \\
			\hline

			{\anita (this paper)}
			& \cellcolor{bgcolor}
			$O\left(n\right)$
			&\cellcolor{bgcolor}  \makecell{$O\Big(n\big(1+\log\frac{1}{\epsilon\sqrt{n}}\big)\Big)$} 
			& \cellcolor{bgcolor} $O\left(n\log\sqrt{n}\right)$ 
			&  \cellcolor{bgcolor} $O\Big(\sqrt{\frac{nL}{\epsilon}}\Big)$\\
			\hline	
			
			\makecell{Lower bound \\ \citep{woodworth2016tight}}
			& $\Omega\left(n\right)$
			& $\Omega\left(n\right)$
			& \makecell{$\Omega\Big(n\sqrt{\frac{L}{\epsilon n}}\Big)$}
			& $\Omega\left(\sqrt{\frac{nL}{\epsilon}}\right)$ \\
			\hline
			
		\end{tabular}
	\end{spacing}
	\vspace{1mm}
	\begin{minipage}{0.99\textwidth}
		\scriptsize
		\centering
	 	{\bf Remark:}  \anita achieves the optimal result $O(n)$ for large-scale problems (large $n$) or moderate error (not too small $\epsilon$). 
		It should be pointed out that all parameter settings of \anita (i.e., $\{p_t\}$, $\{\theta_t\}$, $\{\eta_t\}$, and $\{\alpha_t\}$ in Algorithm \ref{alg:anita}) do not require the value of $\epsilon$ in advance. The convergence rate of \anita will automatically switch to different results listed in Table~\ref{tab:2}.
	\end{minipage}
	\vspace{-10mm}
\end{table}

\section{Our Contributions}

In this paper, we mainly focus on further improving the convergence result in order to close the gap between the upper and lower bound.
We propose a novel loopless accelerated variance-reduced gradient method, called \anita (Algorithm~\ref{alg:anita}), for solving both general convex and strongly convex finite-sum problems given in the form of  \eqref{eq:prob}.
The proposed \anita takes an important step towards the ultimate limit of accelerated methods and can achieve the optimal convergence rates.
Tables~\ref{tab:1} and \ref{tab:2} summarize the convergence results of previous algorithms and \anita.

Now, we highlight the following results achieved by \anita: 

$\bullet$ For general convex problems, \anita obtains the rate $O\big(n\min\big\{1+\log\frac{1}{\epsilon\sqrt{n}},\ \log\sqrt{n}\big\} + \sqrt{\frac{nL}{\epsilon}} \big)$ for finding an $\epsilon$-approximate solution of problem \eqref{eq:prob}, which improves previous best result $O\big(n\min\{\log\frac{1}{\epsilon},\ \log n\} + \sqrt{\frac{nL}{\epsilon}}\big)$ given by Varag~\citep{lan2019unified} (see the `general convex' column of Table~\ref{tab:1}).
Moreover, for a very wide range of $\epsilon$, i.e., $\epsilon \in (0,\frac{L}{n\log^2\sqrt{n}}]\cup [\frac{1}{\sqrt{n}},+\infty)$,
or the number of data samples $n\in (0, \frac{L}{\epsilon\log^2\sqrt{n}}] \cup [\frac{1}{\epsilon^2}, +\infty)$, \anita can exactly achieve the optimal convergence result $O\big(n+\sqrt{\frac{nL}{\epsilon}}\big)$ matching the lower bound $\Omega\big(n+\sqrt{\frac{nL}{\epsilon}}\big)$ provided by \citet{woodworth2016tight} (see Table~\ref{tab:1} and its Footnote~2). 

$\bullet$ In particular, we would like to point out that none of previous algorithms with/without acceleration can obtain the optimal result $O(n)$ for finite-sum problems~\eqref{eq:prob} where the number of data samples is very large or the convergence error is not very small, \anita is the first algorithm that achieves the optimal result $O(n)$ for these typical machine learning problems (see the second column of Table~\ref{tab:2} and its Remark).

$\bullet$ We also note that \anita is the first loopless direct accelerated stochastic algorithm for solving general convex finite-sum problems, while previous accelerated stochastic algorithms use indirect approaches (RPDG, Catalyst, Katyusha) and/or use inconvenient double-loop algorithmic structures (Katyusha\textsuperscript{ns}, Varag) (see Table~\ref{tab:1}). Moreover, by exploiting the loopless structure of \anita, we provide a novel \emph{dynamic multi-stage convergence analysis} which is the key technical part for improving previous results to the optimal rates. 

$\bullet$ For strongly convex finite-sum problems (i.e., under strong convexity Assumption~\ref{asp:strong}), we also prove that \anita achieves the optimal convergence rate $O\big((n+\sqrt{\frac{nL}{\mu}})\log\frac{1}{\epsilon}\big)$ matching the lower bound  $\Omega\big((n+\sqrt{\frac{nL}{\mu}})\log\frac{1}{\epsilon}\big)$ provided by \citet{lan2015optimal} (see Table~\ref{tab:1}).

$\bullet$ Finally, the numerical experiments show that \anita converges faster than the previous state-of-the-art Varag \citep{lan2019unified}, validating our theoretical results and confirming the practical superiority of \anita.

\subsection{\anita algorithm}
\label{sec:alg}
In this section, we describe the simple \anita method in Algorithm~\ref{alg:anita}.
\anita uses the SVRG gradient estimator \eqref{eq:svrg} (see Line \ref{line:grad} of Algorithm~\ref{alg:anita}) and two interpolation steps (momentum) (see Line \ref{line:uxt} and Line \ref{line:bxtn} of Algorithm~\ref{alg:anita}).
Line \ref{line:update} of Algorithm~\ref{alg:anita} is a gradient update step. 

Although previous accelerated stochastic algorithms such as Katyusha/Katyusha\textsuperscript{ns} \citep{allen2017katyusha} and Varag \citep{lan2019unified} also adopt the SVRG gradient estimator combined with momentum steps, \anita enjoys a simpler loopless algorithmic structure.
Note that the previous loopless SVRG/Katyusha algorithms provided in \cite{kovalev2019don} only solve the strongly convex case ($\mu>0$).
Here, our loopless algorithm \anita can deal with both general convex ($\mu=0$) and strongly convex ($\mu>0$) problems, and the \anita algorithm itself is also different and more concise than the loopless algorithms in \cite{kovalev2019don}.
Moreover, for general convex problems ($\mu=0$), \anita provides a new state-of-the-art convergence result which improves all previous results.

\begin{algorithm}[t]
	\caption{\anita (AN optImal loopless acceleraTed vAriance-reduced gradient method)}
	\label{alg:anita}
	\begin{algorithmic}[1]
		\REQUIRE ~
		initial point $x_0$,  parameters $\{p_t\}$, $\{\theta_t\}$, $\{\eta_t\}$, $\{\alpha_t\}$
		\STATE $w_0=\bar{x}_0=\underline{x}_0=x_0$
		\FOR{$t=0,1,2, \ldots, T-1$}
		\STATE $\uxt = \theta_t \xt + (1-\theta_t)\wt$ \label{line:uxt}
		\STATE Randomly pick $i\in\{1,2,\dots,n\}$ \label{line:i}
		\STATE $\tnabla = \nabla f_i(\uxt) - \nabla f_i(\wt) + \nabla f(\wt)$  \label{line:grad}
		\STATE $\xtn = \frac{1}{1+\mu\eta_t} (\xt + \mu \eta_t \uxt)-\frac{\eta_t}{\alpha_t} \tnabla$  \label{line:update}
		\STATE $\bxtn = \theta_t \xtn + (1-\theta_t)\wt$ \label{line:bxtn}
		\STATE $\wtn= \begin{cases}
		\bxtn &\text {with probability } p_t\\
		\wt &\text {with probability } 1-p_t
		\end{cases}$ \label{line:prob}  
		\ENDFOR
		
		\ENSURE $w_T$
	\end{algorithmic}		
\end{algorithm}		

In each iteration $t$, the stochastic gradient estimator $\tnabla$ of \anita (Line \ref{line:grad} of Algorithm~\ref{alg:anita}) uses the gradient information of only one randomly sampled function $f_i$. Note that for the last term $\nabla f(\wt)$, it reuses previous $\nabla f(w_{t-1})$ with probability $1-p_{t-1}$ or needs to compute the full gradient $\nabla f(\bxt)$ with probability $p_{t-1}$ (see Line \ref{line:prob} of Algorithm~\ref{alg:anita}).
Thus we know that \anita uses $(n+2)p_{t-1}+2(1-p_{t-1})$ stochastic gradients in expectation for iteration $t$.
In particular, if $p_t\equiv \frac{1}{n}$, then \anita only uses constant stochastic gradients for each iteration which maintains the same computational cost as SGD.
The snapshot point $w_t$ is updated in the last Line \ref{line:prob} of Algorithm~\ref{alg:anita}, it is a probabilistic step which is the key part for removing double-loop structures to obtain a simple loopless algorithm, similar to \citep{kovalev2019don,li2021page}.
However, our key technical part is that we propose a new dynamic multi-stage convergence analysis which uses a dynamic control of the probability $\{p_t\}$ in Line \ref{line:prob}, unlike directly fixing it to a constant $p_t\equiv p$ as in \citep{kovalev2019don,li2021page}.
To the best of our knowledge, this is the first time that a loopless algorithm uses a dynamic change of $\{p_t\}$.

\section{Preliminaries}
\label{sec:pre}

\topic{Notation} 
Let $[n]$ denote the set $\{1,2,\cdots,n\}$ and $\n{\cdot}$ denote the Euclidean norm for a vector and the spectral norm for a matrix.
Let $\inner{u}{v}$ denote the inner product of two vectors $u$ and $v$.
We use $O(\cdot)$ and $\Omega(\cdot)$ to hide the absolute constant.
We will write $x^* := \arg\min_{x\in \R^d} f(x)$.

For convex problems, one typically uses the function value gap as the 
convergence criterion.
\begin{definition}\label{def:1}
	A point $\hx$ is called an $\epsilon$-approximate solution for problem \eqref{eq:prob} if 
	$\E[f(\hx)-f(x^*)] \leq \epsilon$.
\end{definition}

To show the convergence results, we assume the following standard smoothness assumption for the component functions $f_i$s in \eqref{eq:prob}.
\begin{assumption}[$L$-smoothness]\label{asp:smooth}
	Functions $f_i:\R^d\to \R$ are convex and $L$-smooth such that
	\begin{equation}\label{eq:smooth}
	\n{\nabla f_i(x) - \nabla f_i(y)}\leq L \n{x-y}
	\end{equation}
	for some $L\geq 0$ and all $i\in[n]$.
\end{assumption}
It is easy to see that $f(x)= \frac{1}{n}\sum_{i=1}^n{f_i(x)}$ is also $L$-smooth under Assumption~\ref{asp:smooth}.

For considering the strongly convex setting, we assume the following Assumption~\ref{asp:strong}. 
\begin{assumption}[$\mu$-strong convexity]\label{asp:strong}
	A function $f:\R^d\to \R$ is $\mu$-strongly convex such that
	\begin{equation}\label{eq:strong}
	f(x) - f(y) - \inner{\nabla f(y)}{x-y} \geq \frac{\mu}{2}\ns{x-y},
	\end{equation}
	for some $\mu \geq 0$. 
\end{assumption}
Note that the strong convexity is only corresponding to the average function $f$ in \eqref{eq:prob}, is not needed for the component functions $f_i$s.

\section{Convergence Results for \anita}
\label{sec:result}

In this section, we present two main convergence theorems of \anita (Algorithm~\ref{alg:anita}) for solving finite-sum problems \eqref{eq:prob}, i.e., Theorem~\ref{thm:1} (general convex setting in Section \ref{sec:general}) and Theorem~\ref{thm:strong} (strongly convex setting in Section \ref{sec:strong}).
Subsequently, we formulate two Corollaries~\ref{cor:1}--\ref{cor:strong} from Theorems~\ref{thm:1}--\ref{thm:strong} for providing the detailed convergence results. 
The detailed proofs for Theorems~\ref{thm:1}--\ref{thm:strong} and Corollaries~\ref{cor:1}--\ref{cor:strong} are deferred to Appendix \ref{sec:app-proofs}.

\subsection{General convex setting}
\label{sec:general}

In this section, we provide the main convergence theorem of \anita for general convex problems and then obtain a corollary for providing the detailed convergence result. 
Note that if we fix the probability $p_t$ in Line \ref{line:prob} of Algorithm~\ref{alg:anita} to a constant $p$, then the update of $w_t$ follows from a geometric distribution $\text{Geom}(p)$.
For a geometric distribution $N\sim \text{Geom}(p)$, i.e., $N=k$ with probability $(1-p)^kp$ for $k=0,1,2,\ldots$ (after $k$ failures until the first success), we know that $\E[N] = \frac{1-p}{p}$. 
In the \emph{first stage} of \anita, we indeed use constant probability $p_t\equiv p= \frac{1}{n+1}$. Let $t_1$ be the first time such that $w$ changes to $\bar x$, i.e., $w_{t_1+1}=\bar x_{t_1+1}$ and $w_{t_1}=w_{t_1-1}=\cdots=w_0$.
Thus $t_1\sim \text{Geom}(p)$ and $\E[t_1] =n$.
Note that this first stage where we fix $p_t\equiv p$  is similar to loopless SVRG \citep{kovalev2019don}, SCSG \citep{lei2016less} and PAGE \citep{li2021page}. 
One can also derandomize the special case of constant probability $p$ in this first stage to a deterministic double-loop with loop length $\frac{1-p}{p}$ algorithms like the original SVRG \citep{johnson2013accelerating} and SARAH \citep{nguyen2017sarah}.
The difference is that our \anita will use a dynamic change of $p_t$ after the first stage, while previous algorithms always keep fixing the probability $p_t\equiv p$.

\begin{theorem}[General convex case]\label{thm:1}
	Suppose that Assumption~\ref{asp:smooth} holds. 
	For $0\leq t\leq t_1$, let $p_t\equiv \frac{1}{n+1}$, $\theta_t \equiv 1-\frac{1}{2\sqrt{n}}$, $\etat \leq  \frac{1}{L(1+1/(1-\theta_t))}$ and $\alpha_t =\theta_t$.
	For $t> t_1$, let $p_t=\max\{\frac{4}{t-t_1+3\sqrt{n}}, \frac{4}{n+3}\}$, $\theta_t =\frac{2}{p_t(t-t_1+3\sqrt{n})}$, $\etat \leq  \frac{1}{3L}$ and $\alpha_t =\theta_t$. 
	Then the following equation holds for \anita (Algorithm~\ref{alg:anita}) for any iteration $t>t_1+1$:
	\begin{align*}
	\E[f(w_t)-f(x^*)] 
	\leq  	\frac{32\ns{x_0-x^*}}{\eta_{t-1}p_{t-1}(t-t_1+3\sqrt{n})^2}.
	\end{align*}
\end{theorem}

According to Theorem~\ref{thm:1}, we can obtain a detailed convergence result in the following Corollary~\ref{cor:1}.

\begin{corollary}[General convex case]\label{cor:1}
	Suppose that Assumption~\ref{asp:smooth} holds. 
	Choose the parameters $\{p_t\}$, $\{\theta_t\}$, $\{\eta_t\}$, $\{\alpha_t\}$ as stated in Theorem~\ref{thm:1}.
	Then \anita (Algorithm~\ref{alg:anita}) can find an $\epsilon$-approximate solution for problem \eqref{eq:prob} such that
	$$\E[f(w_T)-f(x^*)]\leq \epsilon$$ 
	within $T$ iterations, where 
	$$T\leq \begin{cases}
	2n                  & \text{if }~ \epsilon \geq O(\frac{1}{n})\\
	n+\sqrt{\frac{24(n+3)L\ns{x_0-x^*}}{\epsilon}} &\text{if }~\epsilon < O(\frac{1}{n})
	\end{cases},$$ 
	and the number of stochastic gradient computations can be bounded by
	\begin{align*}
	\#\mathrm{grad} =O\left(n \min\Big\{1+\log\frac{1}{\epsilon\sqrt{n}},\ \log\sqrt{n}\Big\} + \sqrt{\frac{nL}{\epsilon}} \right).
	\end{align*}
\end{corollary}

\topic{Remark}
From the choice of probability $\{p_t\}$ in Theorem~\ref{thm:1}, we know that there are three stages of \anita: i) the first stage $p_t\equiv\frac{1}{n+1}$ for $0\leq t\leq t_1$;  ii) the second stage $p_t=\frac{4}{t-t_1+3\sqrt{n}}$ for $t_1<t \leq t_1+n +3-3\sqrt{n}$; iii) the third stage $p_t\equiv \frac{4}{n+3}$ for $t>t_1+n +3-3\sqrt{n}$.
This novel multi-stage convergence analysis is key part for the improvement of \anita.
Roughly speaking, the number of stochastic gradient computations in the first stage is $\#\mathrm{grad} =O(n)$, in the second stage is $\#\mathrm{grad} =O\big(n\min\big\{\log\frac{1}{\epsilon\sqrt{n}}, \log\sqrt{n}\big\}\big)$, and in the third stage is $\#\mathrm{grad} =O\big(\sqrt{\frac{nL}{\epsilon}}\big)$.
The detailed proofs of Theorem~\ref{thm:1} and Corollary~\ref{cor:1} are deferred to Appendix \ref{sec:app-1}.
Note that the guarantee of \anita is the \emph{last iterate} convergence  unlike previous \emph{average iterates} convergence. 
Also note that all parameter settings $\{p_t\}$, $\{\theta_t\}$, $\{\eta_t\}$, $\{\alpha_t\}$ of \anita in Theorem~\ref{thm:1} do not require the value of $\epsilon$ in advance. 
The convergence rate of \anita will automatically switch to different results as stated in Table~\ref{tab:2}.

\subsection{Strongly convex setting}
\label{sec:strong}

In this section, we provide the main convergence theorem of \anita for strongly convex problems ($\mu>0$ in Assumption~\ref{asp:strong}) and then obtain a corollary for providing the detailed convergence result. 

\begin{theorem}[Strongly convex case]\label{thm:strong}
	Suppose that Assumptions \ref{asp:smooth} and \ref{asp:strong} hold. 
	For any $t\geq0$, let $p_t\equiv p$, $\theta_t \equiv \theta=\frac{1}{2}\min\{1,\sqrt{\frac{\mu}{pL}}\}$, 
	$\etat \leq  \frac{1}{L\theta_t(1+1/(1-\theta_t))}$ and $\alpha_t =1+\mu\eta_t$. 
	Then the following equation holds for \anita (Algorithm~\ref{alg:anita}) for any iteration $t\geq 0$:  
	\begin{align}
	\E[\Phi_t] \leq \Big(1-\frac{4p\theta}{5}\Big)^t \Phi_0, \label{eq:linear}
	\end{align}
	where $\Phi_t:=f(w_t)-f(x^*)+\frac{(1+\mu\eta)p\theta}{2\eta}\ns{x_t-x^*}$.
\end{theorem}

Similarly, according to Theorem~\ref{thm:strong}, we can obtain a detailed convergence result in the following Corollary~\ref{cor:strong}.

\begin{corollary}[Strongly convex case]\label{cor:strong}
	Suppose that Assumptions \ref{asp:smooth} and \ref{asp:strong} hold. 
	Choose the parameters $\{p_t\}$, $\{\theta_t\}$, $\{\eta_t\}$, $\{\alpha_t\}$ as stated in Theorem~\ref{thm:strong}.
	Then \anita (Algorithm~\ref{alg:anita}) can find an $\epsilon$-approximate solution for problem \eqref{eq:prob} such that
	$$\E[f(w_T)-f(x^*)]\leq \epsilon$$ 
	within $T$ iterations, where 
	$$T\leq \frac{5}{4p\theta}\log\frac{\Phi_0}{\epsilon}.$$ 
	Moreover, by choosing $p=\frac{1}{n}$ and recalling that $\theta=\frac{1}{2}\min\{1,\sqrt{\frac{\mu}{pL}}\}$, the number of stochastic gradient computations can be bounded by
	\begin{align*}
	\#\mathrm{grad} =O\left(\max\Bigg\{n, \sqrt{\frac{nL}{\mu}}\Bigg\}\log\frac{1}{\epsilon}\right).
	\end{align*}
\end{corollary}

\topic{Remark}
In this strongly convex case, the parameter setting of \anita in Theorem~\ref{thm:strong} is simpler than the general convex case in Theorem~\ref{thm:1}. 
Here, the choice of probability $\{p_t\}$ can be fixed to a constant $p$ and $\{\theta_t\}$ also can be chosen as a constant $\theta$.
Then according to Theorem~\ref{thm:strong}, we know that $\{\eta_t\}$ and $\{\alpha_t\}$ also reduce to constant values.
Thus there is only one stage in this strongly convex case rather than three stages in previous general convex case. 
Also here the function value decreases in an exponential rate, i.e., $\E[\Phi_t] \leq \big(1-\frac{4p\theta}{5}\big)^t \Phi_0$ (see \eqref{eq:linear} in Theorem~\ref{thm:strong}).
It is easy to see that the number of iterations $T$ can be bounded by $O(\cdot \log\frac{1}{\epsilon})$ for finding an $\epsilon$-approximate solution $\E[f(w_T)-f(x^*)]\leq \epsilon$.
Then, by choosing $p=\frac{1}{n}$ (thus each iteration only computes constant stochastic gradients in expectation), the number of total stochastic gradient computations can be bounded by
$\#\mathrm{grad} =O\big(\max\big\{n, \sqrt{\frac{nL}{\mu}}\big\}\log\frac{1}{\epsilon}\big)$.
This convergence result is optimal which matches the lower bound  $\Omega\Big(\big(n+\sqrt{\frac{nL}{\mu}}\big)\log\frac{1}{\epsilon}\Big)$ given by \citet{lan2015optimal} (see Table~\ref{tab:1}).
The detailed proofs of Theorem~\ref{thm:strong} and Corollary~\ref{cor:strong} are deferred to Appendix \ref{sec:app-strong}.
Note that all parameter settings $\{p_t\}$, $\{\theta_t\}$, $\{\eta_t\}$, $\{\alpha_t\}$ of \anita in Theorem~\ref{thm:strong} also do not require the value of $\epsilon$ in advance.

\section{Experiments}
In this section, we present the numerical experiments of \anita (Algorithm~\ref{alg:anita}) compared with previous state-of-the-art Varag~\citep{lan2019unified}. We also present the standard gradient descent (GD) as a benchmark for demonstrating the performance of these algorithms.
The theoretical convergence results of these algorithms can be found in Table~\ref{tab:1}.

In the experiments, we consider the following logistic  regression problem:
\begin{equation}\label{eq:prob-exp}
\min_{x\in \R^d}   f(x):= \frac{1}{n}\sum_{i=1}^n \log\big(1+\exp(-b_ia_i^\top x)\big),
\end{equation}
where $\{a_i,b_i\}_{i=1}^n\in \R^{d}\times \{\pm 1\}$ are data samples.
All datasets used in our experiments are downloaded from LIBSVM \citep{chang2011libsvm}.
We also point out that we directly use the parameter settings according to the theoretical convergence theorems or corollaries of these algorithms, i.e., we do not tune any hyperparameters.
Note that for the logistic function in \eqref{eq:prob-exp}, one can precompute the smoothness parameter $L$ satisfying Assumption~\ref{asp:smooth}, i.e., $L\leq 1/4$ if the data samples are normalized.
Given the parameter $L$, we are ready to set all other hyperparameters for GD (see Corollary 2.1.2 in \citealp{nesterov2014introductory}), for Varag (see Theorem 1 in \citealp{lan2019unified}) and for \anita (see our Theorem~\ref{thm:1}).
Note that all of these three algorithms only require $L$ for setting their (hyper)parameters.

In the following Figure~\ref{fig:1}, the $x$-axis and $y$-axis represent the number of data passes (i.e., we compute $n$ stochastic gradients for each data pass) and the training loss, respectively. The numerical results presented in Figure~\ref{fig:1} are conducted on different datasets. Each plot corresponds to one dataset (six datasets in total). 
The experimental results show that \anita indeed converges faster than Varag \citep{lan2019unified} in the earlier stage (moderate convergence error), validating our theoretical results (see the second column of Table~\ref{tab:2} and its Remark).
More importantly, \anita is the first accelerated algorithm which can obtain the optimal convergence result $O(n)$ in this range.
Besides, \anita also enjoys a simpler loopless algorithmic structure while Varag uses an double-loop structure.	

\begin{figure}[h]
	\centering
	\begin{minipage}{0.34\textwidth}
		\centering
		\includegraphics[width=\linewidth]{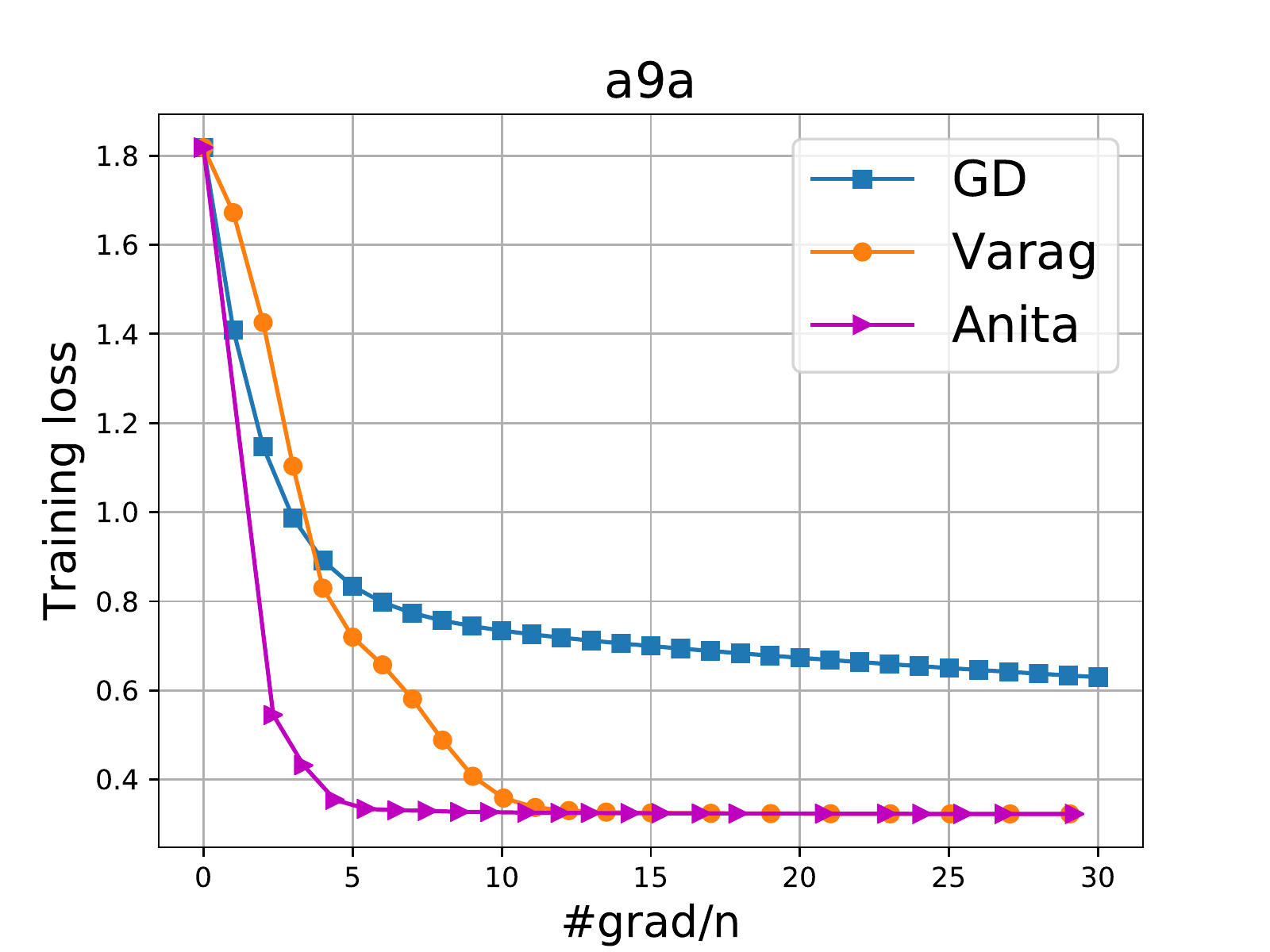}
	\end{minipage}%
	\begin{minipage}{0.34\textwidth}
		\centering
		\includegraphics[width=\linewidth]{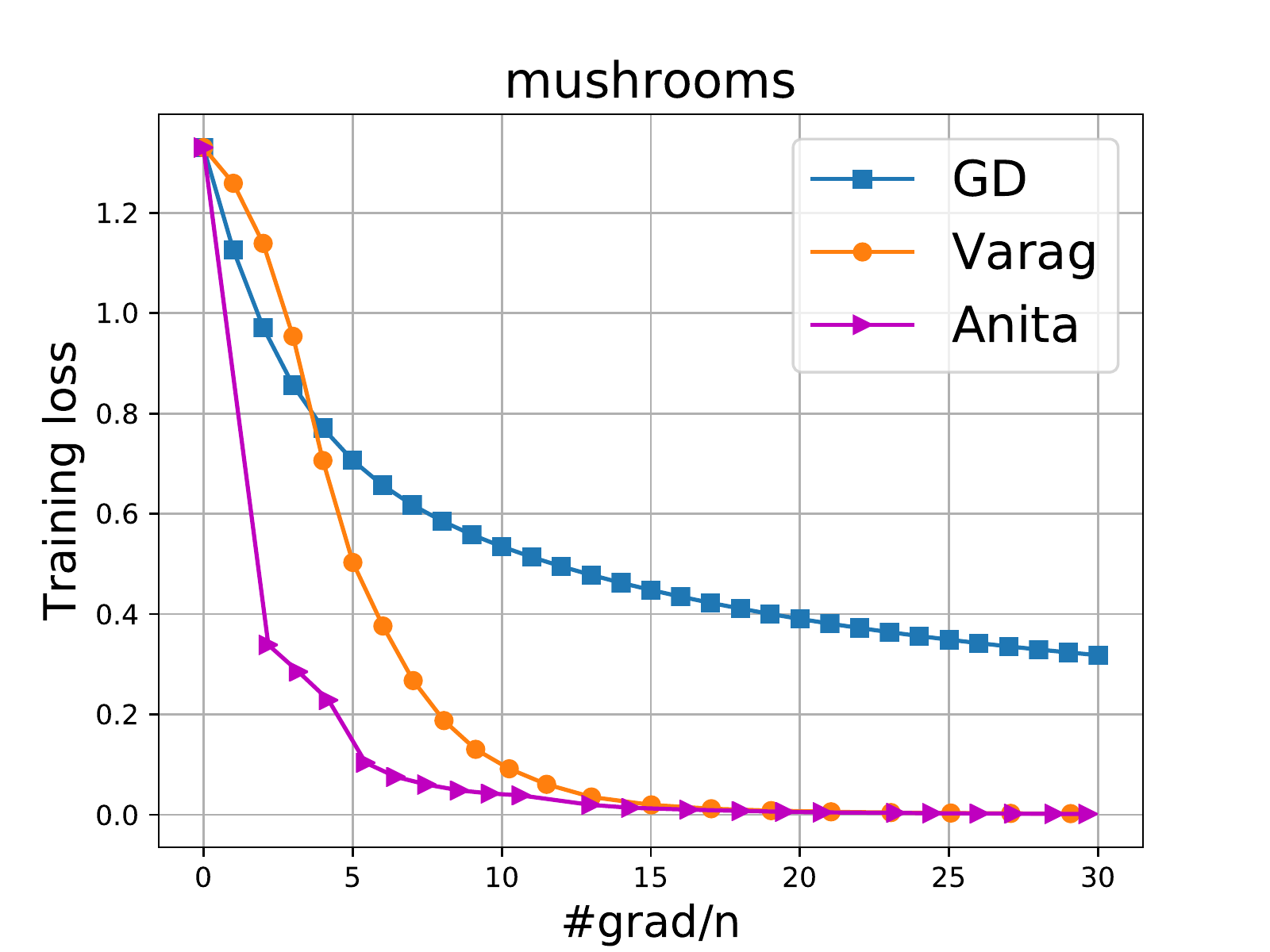}
	\end{minipage}%
	\begin{minipage}{0.34\textwidth}
		\centering
		\includegraphics[width=\linewidth]{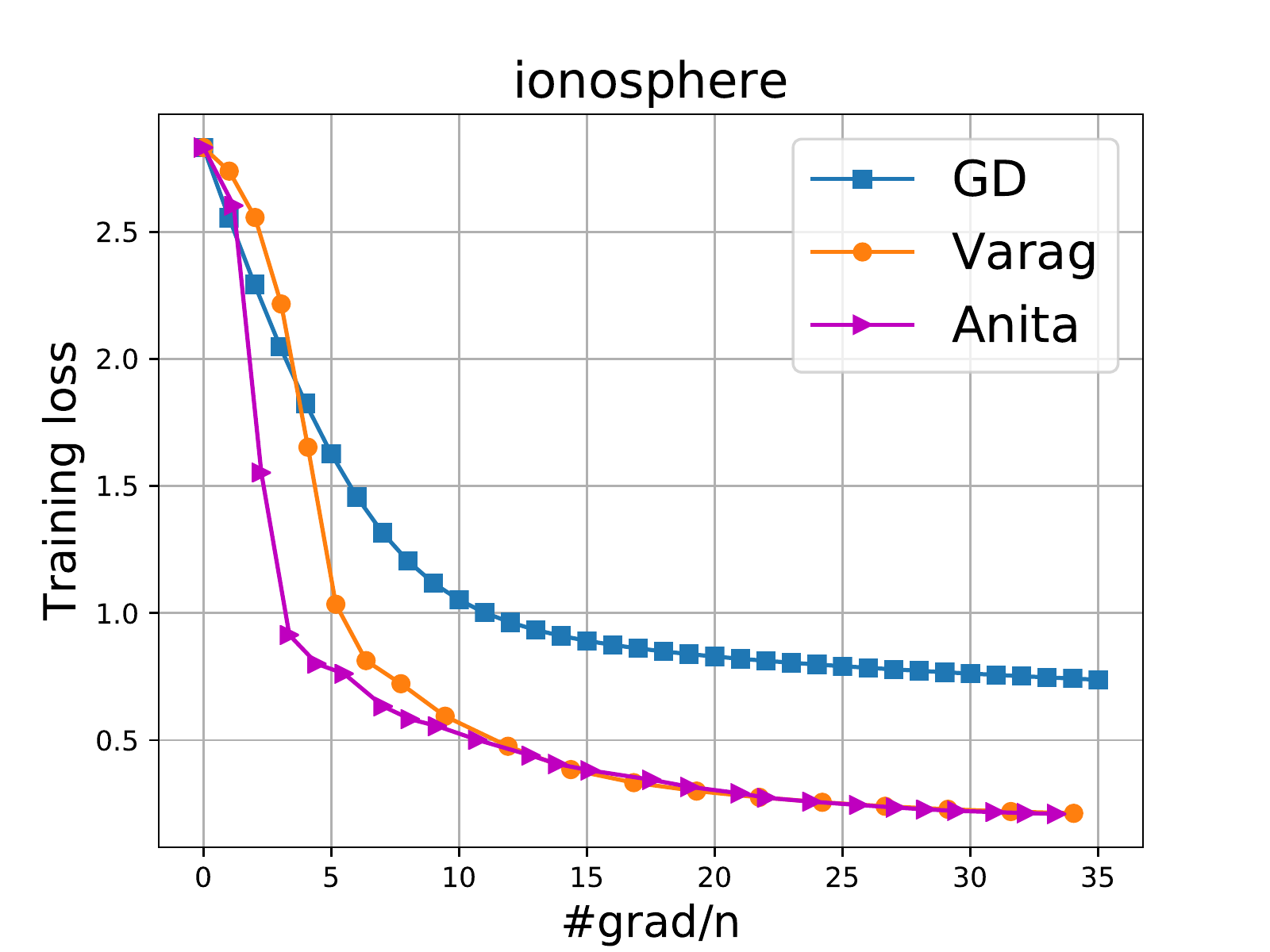}
	\end{minipage}\\
	\vspace{2mm}
	\begin{minipage}{0.34\textwidth}
		\centering
		\includegraphics[width=\linewidth]{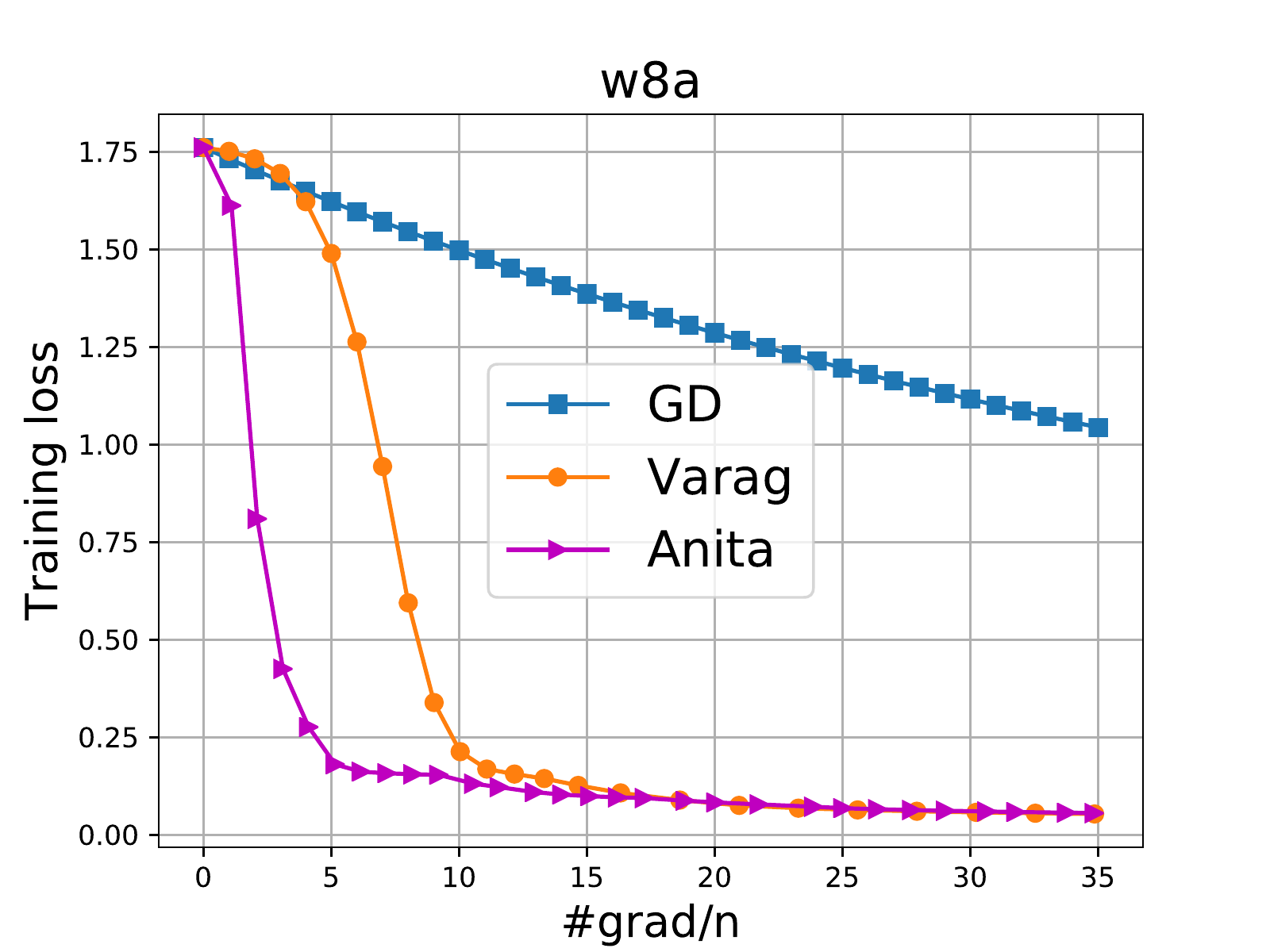}
	\end{minipage}%
	\begin{minipage}{0.34\textwidth}
		\centering
		\includegraphics[width=\linewidth]{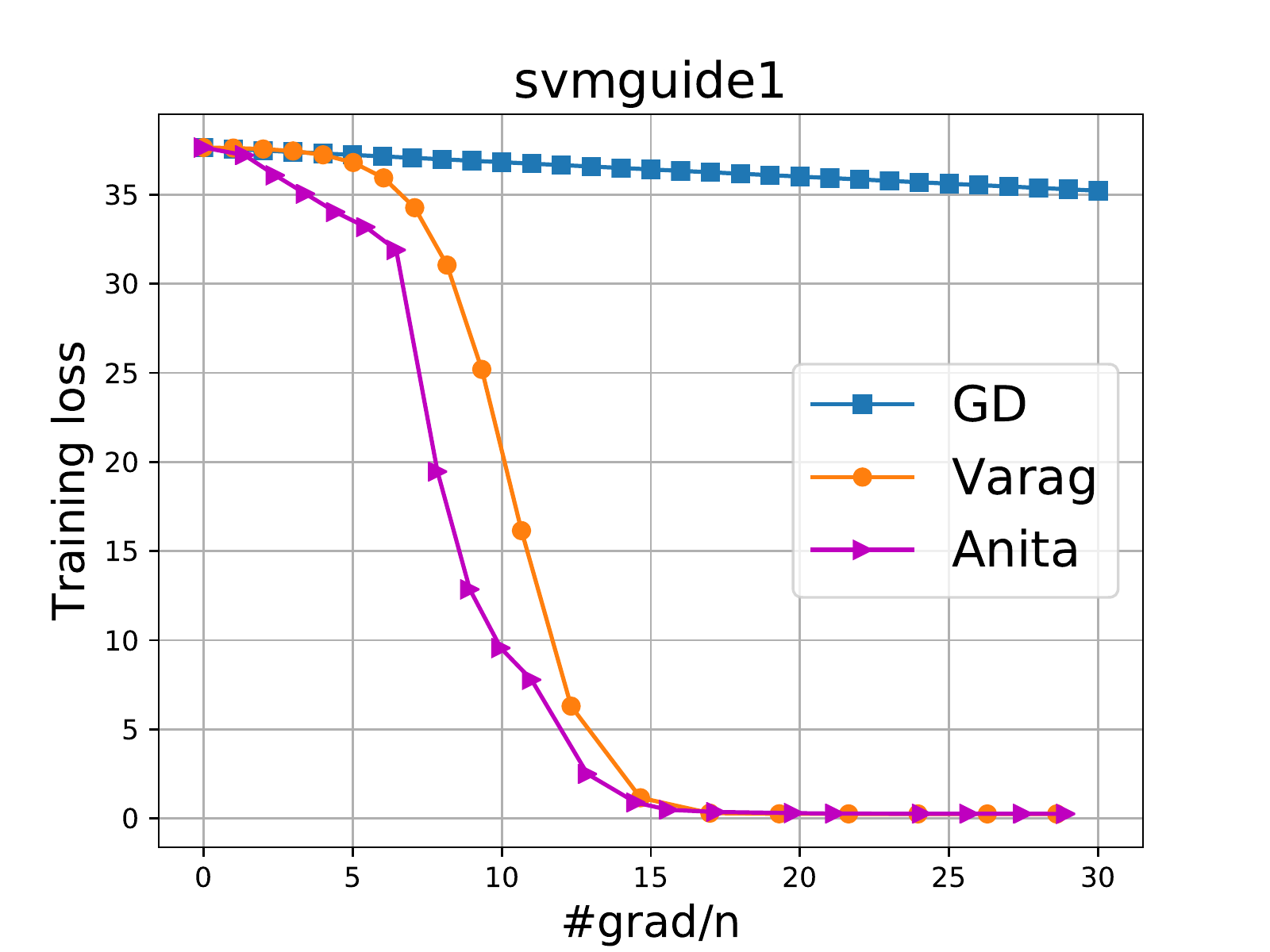}
	\end{minipage}%
	\begin{minipage}{0.34\textwidth}
		\centering
		\includegraphics[width=\linewidth]{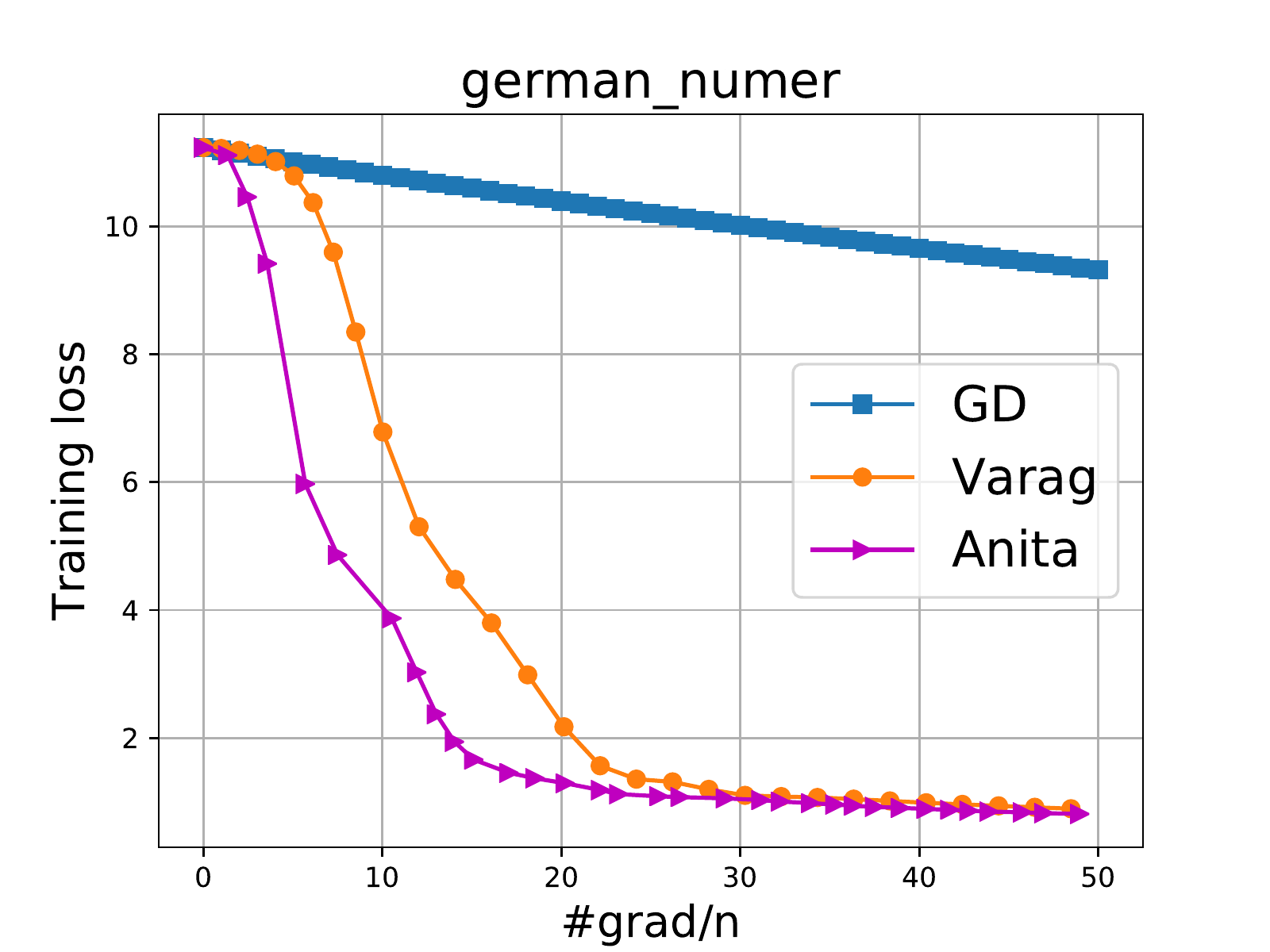}
	\end{minipage}\\
	\caption{The convergence performance of GD, Varag and \anita under different datasets.}
	\label{fig:1}
\end{figure}	

\section{Conclusion}
In this paper, we propose a simple accelerated variance-reduced gradient method \anita, for solving both general convex and strongly convex finite-sum problems. 
The proposed \anita takes an important step towards the ultimate limit of accelerated methods to close the gap between the upper and lower bound.
In particular, \anita achieves the first optimal convergence rate $O(n)$ matching the lower bound $\Omega(n)$ for large-scale general convex problems.
Besides, it also achieves the optimal convergence rate $O\big((n+\sqrt{\frac{nL}{\mu}})\log\frac{1}{\epsilon}\big)$ matching the lower bound  $\Omega\big((n+\sqrt{\frac{nL}{\mu}})\log\frac{1}{\epsilon}\big)$ for strongly convex problems.
Moreover, we provide a novel dynamic multi-stage convergence analysis utilizing the simpler loopless algorithmic structure, which is the key technical part for improving previous results to the optimal rates. 
The numerical experiments validate our theoretical results and confirm the practical superiority of \anita.
Our new theoretical rates and convergence analysis can also lead to key improvements for many other related problems.
For instance, \citet{li2021canita} obtain the first accelerated result, substantially improving previous state-of-the-art results, by applying \anita to the distributed optimization problems with compressed communication.


\bibliographystyle{plainnat}
\bibliography{anita}

\newpage

\appendix
\section{Missing Proofs}
\label{sec:app-proofs}
Now, we provide the detailed proofs of main convergence theorems and corollaries of \anita for both general convex case (Theorem~\ref{thm:1} and Corollary~\ref{cor:1}) and strongly convex case (Theorem~\ref{thm:strong} and Corollary~\ref{cor:strong}).

Before proving these theorems and corollaries, we first recall some basic properties for smooth convex functions (e.g., \citealp{nesterov2014introductory,lan2019unified}) and some basic facts for the geometric distribution (e.g., \citealp{lei2017non}).

\begin{lemma}[Lemma 1 in \citealp{lan2019unified}] \label{lem:smoothfunc}
	If $f: X \to \R$ has $L$-Lipschitz continuous gradients ($L$-smooth), then  we have 
	\begin{align}
	\frac{1}{2L} \|\nabla f(x) - \nabla f(z)\|^2 \leq
	f(x) - f(z) - \langle \nabla f(z), x-z \rangle,  \qquad \forall x, z \in X. \label{eq:l_varag}
	\end{align}
\end{lemma}

\begin{lemma}[Lemma A.2 in \citealp{lei2017non}] \label{lem:geom}
	Let $N\sim \mathrm{Geom}(p)$.  Then for any sequence $D_0, D_1,\ldots$ with $\E|D_N|<\infty$, we have 
	\begin{align}
	&\E[N] = \frac{1-p}{p},  \label{eq:factexp}\\
	&\E[D_N-D_{N+1}] = \frac{p}{1-p}\big(D_0-\E[D_N]\big),   \label{eq:factn1} \\
	&\E[D_N] = pD_0 + (1-p)\E[D_{N+1}].  \label{eq:factn}
	\end{align}
\end{lemma}

Now, we provide some important technical lemmas which are useful for proving the main convergence theorems of \anita.
Concretely, Lemma~\ref{lem:var} provides some ways to upper bound the variance of the gradient estimator in \anita.
Lemma~\ref{lem:func_relation} describes the change of function value after a gradient update step in \anita.

\begin{lemma}\label{lem:var}
	Suppose that Assumption~\ref{asp:smooth} holds.  
	The gradient estimator 
	\begin{align}\label{eq:tnabla}
	\tnabla = \nabla f_i(\uxt) - \nabla f_i(\wt) + \nabla f(\wt)
	\end{align}
	is defined in Line \ref{line:grad} of Algorithm~\ref{alg:anita}, then conditional on the past, we have 
	\begin{align}
	&\E[\tnabla] = \nabla f(\uxt), \label{eq:expvar}\\
	&\E[\ns{\tnabla-\nabla f(\uxt)}] \leq L^2 \ns{\uxt-\wt}, \label{eq:var1}\\
	&\E[\ns{\tnabla-\nabla f(\uxt)}] \leq 2L 
	\big(f(\wt) - f(\uxt)  - \inner{\nabla f(\uxt)}{\wt-\uxt}\big). \label{eq:var2}
	\end{align}
\end{lemma}
\begin{proofof}{Lemma~\ref{lem:var}}
	For \eqref{eq:expvar}, it is easy to see that (note that the expectation is taken over the random choice of $i$ in iteration $t$ (see Line \ref{line:i} of Algorithm~\ref{alg:anita})) 
	\begin{align*}
	\E[\tnabla] 
	&\overset{\eqref{eq:tnabla}}{=}\E[\nabla f_i(\uxt) - \nabla f_i(\wt) + \nabla f(\wt)] \notag\\
	&=\nabla f(\uxt) - \nabla f(\wt) + \nabla f(\wt)
	=\nabla f(\uxt).
	\end{align*}
	Then, for \eqref{eq:var1}, we obtain it from Assumption~\ref{asp:smooth} as follows:
	\begin{align}
	\E[\ns{\tnabla-\nabla f(\uxt)}] 
	&\overset{\eqref{eq:tnabla}}{=}\E[\ns{\nabla f_i(\uxt) - \nabla f_i(\wt) + \nabla f(\wt) - \nabla f(\uxt)}] \notag \\
	&\leq \E[\ns{\nabla f_i(\uxt) - \nabla f_i(\wt)}] \label{eq:usefactvar} \\
	&\leq L^2\ns{\uxt - \wt}, \label{eq:usesmooth}
	\end{align}
	where \eqref{eq:usefactvar} follows from the fact that $\E[\ns{x-\E x}]\leq \E[\ns{x}]$ for any random variable $x$, and \eqref{eq:usesmooth} follows from Assumption~\ref{asp:smooth}, i.e., the $L$-Lipschitz continuous gradients $\n{\nabla f_i(x) - \nabla f_i(y)}\leq L \n{x-y}$.
	
	Now, for the last one \eqref{eq:var2},  we obtain it from \eqref{eq:usefactvar} and Assumption~\ref{asp:smooth} as follows:
	\begin{align}
	\E[\ns{\tnabla-\nabla f(\uxt)}] 
	& \overset{\eqref{eq:usefactvar}}{\leq}  \E[\ns{\nabla f_i(\uxt) - \nabla f_i(\wt)}] \notag\\
	&\leq \E\big[2L 
	\big(f_i(\wt) - f_i(\uxt)  - \inner{\nabla f_i(\uxt)}{\wt-\uxt}\big)\big]  \label{eq:use-lvarag} \\
	&=2L \big(f(\wt) - f(\uxt)  - \inner{\nabla f(\uxt)}{\wt-\uxt}\big), \notag
	\end{align}
	where \eqref{eq:use-lvarag} uses Lemma~\ref{lem:smoothfunc} with $x$ and $z$ replaced by $\wt$ and $\uxt$, and $f$ replaced by $f_i$ since  $f_i$ has $L$-Lipschitz continuous gradients according to Assumption~\ref{asp:smooth}.
\end{proofof}

\begin{lemma}\label{lem:func_relation}
	Suppose that Assumptions \ref{asp:smooth} and \ref{asp:strong} hold. 
	Let stepsize $\etat \leq  \frac{\alpha_t}{L(1+\mu\eta_t)\theta_t(1+1/(1-\theta_t))}$,
	then the following equation holds for \anita (Algorithm~\ref{alg:anita}) for any iteration $t\geq 0$:  	
	\begin{align}
	\E[f(w_{t+1})-f(x^*) ]  
	& \leq \E\bigg[(1-p_t\theta_t)\big(f(w_{t})-f(x^*)\big) \notag \\
	&\qquad \qquad 
	+ \frac{p_t\alpha_t\theta_t}{(1+\mu\eta_t)\eta_t}\Big(\frac{1}{2}\ns{x_{t}-x^*} -\frac{1+\mu\eta_t}{2}\ns{x_{t+1}-x^*}\Big) \notag \\
	&\qquad \qquad 
	-\frac{\mu(1+\mu\eta_t-\alpha_t)p_t\theta_t}{2(1+\mu\eta_t)}\ns{\uxt-x^*}
	\bigg]. \label{eq:keyrelation}
	\end{align}
	Note that for the case of $\mu=0$ (general (non-strongly) convex setting), only the smoothness Assumption~\ref{asp:smooth} is required, i.e., the strong convexity Assumption~\ref{asp:strong} is not needed for obtaining \eqref{eq:keyrelation} with $\mu=0$.
\end{lemma}
\begin{proofof}{Lemma~\ref{lem:func_relation}}
	First, in view of $L$-smoothness of $f$ (Assumption~\ref{asp:smooth}), we have
	\begin{align}
	&\E[f(\bxtn)]  \notag\\
	&\leq \E\bigg[f(\uxt) + \inner{\nabla f(\uxt)}{\bxtn - \uxt} + \frac{L}{2} \ns{\bxtn-\uxt}\bigg] \notag\\
	&= \E\bigg[f(\uxt) + \inner{\nabla f(\uxt)}{\theta_t(\xtn-\xt)} +\frac{L\theta_t^2}{2}\ns{\xtn-\xt}\bigg] \label{eq:use-two-m} \\
	&= \E\bigg[f(\uxt) + \inner{\nabla f(\uxt) - \tnabla}{\theta_t(\xtn-\xt)} 
	+\inner{\tnabla}{\theta_t(\xtn-\xt)} +\frac{L\theta_t^2}{2}\ns{\xtn-\xt}\bigg]  \notag\\
	&\leq \E\bigg[f(\uxt) + \frac{\beta_t}{2L}\ns{\nabla f(\uxt) - \tnabla} 
	+\frac{L\theta_t^2}{2\beta_t}\ns{\xtn-\xt} 
	+\inner{\tnabla}{\theta_t(\xtn-\xt)} +\frac{L\theta_t^2}{2}\ns{\xtn-\xt}\bigg]   \label{eq:useyoung} \\
	&= \E\bigg[f(\uxt) + \frac{\beta_t}{2L}\ns{\nabla f(\uxt) - \tnabla} 
	+\frac{L(1+1/\beta_t)\theta_t^2}{2}\ns{\xtn-\xt} 
	+\inner{\tnabla}{\theta_t(x^*-\xt)}\notag\\
	&\qquad \qquad \qquad
	+\inner{\tnabla}{\theta_t(\xtn-x^*)} \bigg]  \notag\\
	& \overset{\eqref{eq:expvar}}{=} \E\bigg[f(\uxt) + \frac{\beta_t}{2L}\ns{\nabla f(\uxt) - \tnabla}  
	+\frac{L(1+1/\beta_t)\theta_t^2}{2}\ns{\xtn-\xt}  
	+ \inner{\nabla f(\uxt)}{\theta_t(x^*-\xt)} \notag\\
	&\qquad \qquad \qquad
	+\inner{\tnabla}{\theta_t(\xtn-x^*)} \bigg] \label{eq:reuse} \\
	&\overset{\eqref{eq:var2}}{\leq} 
	\E\bigg[f(\uxt) + \beta_t\big(f(\wt) - f(\uxt)  - \inner{\nabla f(\uxt)}{\wt-\uxt}\big)
	+\frac{L(1+1/\beta_t)\theta_t^2}{2}\ns{\xtn-\xt} \notag\\
	&\qquad \qquad \qquad
	+ \inner{\nabla f(\uxt)}{\theta_t(x^*-\xt)}
	+\inner{\tnabla}{\theta_t(\xtn-x^*)}\bigg]    \notag\\
	& = \E\bigg[(1-\theta_t)f(\wt) +\theta_t f(\uxt)  - \inner{\nabla f(\uxt)}{(1-\theta_t)(\wt-\uxt)} + \inner{\nabla f(\uxt)}{\theta_t(x^*-\xt)}\notag\\
	&\qquad \qquad \qquad
	+\frac{L(1+1/(1-\theta_t))\theta_t^2}{2}\ns{\xtn-\xt} 
	+\inner{\tnabla}{\theta_t(\xtn-x^*)} \bigg]  \label{eq:use-beta} \\
	& = \E\bigg[(1-\theta_t)f(\wt) +\theta_t \Big(f(\uxt)  + \inner{\nabla f(\uxt)}{x^*-\uxt}\Big)  \notag\\
	&\qquad \qquad \qquad
	+\frac{L(1+1/(1-\theta_t))\theta_t^2}{2}\ns{\xtn-\xt} 
	+\inner{\tnabla}{\theta_t(\xtn-x^*)} \bigg], \label{eq:tmp1} 
	\end{align}
	where \eqref{eq:use-two-m} holds since $\bxtn - \uxt = \theta_t(\xtn-\xt)$ according to the two interpolation steps of \anita (see Line \ref{line:uxt} and Line \ref{line:bxtn} of Algorithm~\ref{alg:anita}), \eqref{eq:useyoung} uses Young's inequality with $\beta_t>0$, \eqref{eq:use-beta} holds by further choosing $\beta_t=1-\theta_t$, \eqref{eq:tmp1} removes $\wt$ and $\xt$ via the interpolation step $\uxt=\theta_t\xt+(1-\theta_t)\wt$ (see Line \ref{line:uxt} of Algorithm~\ref{alg:anita}).
	
	Now, we use the (strong) convexity of $f$ (see Assumption~\ref{asp:strong}) in \eqref{eq:tmp1} to obtain 
	\begin{align}
	\E[f(\bxtn)] 
	& \leq \E\bigg[(1-\theta_t)f(\wt) 
	+\theta_t \Big(f(x^*)-\frac{\mu}{2}\ns{\uxt-x^*}\Big)\notag\\
	&\qquad \qquad \qquad
	+\frac{L(1+1/(1-\theta_t))\theta_t^2}{2}\ns{\xtn-\xt} 
	+\inner{\tnabla}{\theta_t(\xtn-x^*)} \bigg]. \label{eq:tmp2}
	\end{align}
	Then, we deduce the last inner product term in \eqref{eq:tmp2} as follows:
	\begin{align}
	&\E\big[\inner{\tnabla}{\theta_t(\xtn-x^*)}\big] \notag\\
	& = \E\bigg[ \frac{\alpha_t\theta_t}{(1+\mu\eta_t)\eta_t} \innerb{\xt+\mu\eta_t \uxt- (1+\mu\eta_t)\xtn}{\xtn-x^*}  \bigg] \label{eq:use-tnabla}\\
	&= \E\bigg[ \frac{\alpha_t\theta_t}{(1+\mu\eta_t)\eta_t} \Big( \innerb{\xt-\xtn}{\xtn-x^*} +\mu\eta_t \innerb{\uxt-\xtn}{\xtn-x^*}\Big) \bigg] \notag\\
	&=\E\bigg[ \frac{\alpha_t\theta_t}{(1+\mu\eta_t)\eta_t} \Big(\frac{1}{2} \big(\ns{\xt-x^*}-\ns{\xt-\xtn}-\ns{\xtn-x^*} \big)
	\notag\\
	&\qquad \qquad \qquad \qquad \quad
	+\frac{\mu\eta_t}{2} \big(\ns{\uxt-x^*}-\ns{\uxt-\xtn}-\ns{\xtn-x^*}\big)\Big) \bigg] \notag\\
	&\leq \E\bigg[ \frac{\alpha_t\theta_t}{(1+\mu\eta_t)\eta_t} \Big(\frac{1}{2}\ns{\xt-x^*}-\frac{1+\mu\eta_t}{2}\ns{\xtn-x^*} 
	+\frac{\mu\eta_t}{2} \ns{\uxt-x^*}
	-\frac{1}{2}\ns{\xt-\xtn}\Big) \bigg], \label{eq:tmp-inner}
	\end{align}
	where \eqref{eq:use-tnabla} follows from the gradient update step of \anita (see Line \ref{line:update} of Algorithm~\ref{alg:anita}).
	
	Now we plug \eqref{eq:tmp-inner} into \eqref{eq:tmp2} to get
	\begin{align}
	&\E[f(\bxtn)]  \notag\\
	& \leq \E\bigg[(1-\theta_t)f(w_{t})+\theta_t f(x^*) 
	-\frac{\mu(1+\mu\eta_t-\alpha_t)\theta_t}{2(1+\mu\eta_t)}\ns{\uxt-x^*}
	+\frac{L(1+1/(1-\theta_t))\theta_t^2}{2}\ns{\xtn-\xt} \notag \\
	&\qquad \qquad
	+ \frac{\alpha_t\theta_t}{(1+\mu\eta_t)\eta_t}\Big(\frac{1}{2}\ns{x_{t}-x^*} -\frac{1+\mu\eta_t}{2}\ns{x_{t+1}-x^*}\Big)
	-\frac{\alpha_t\theta_t}{2(1+\mu\eta_t)\eta_t}\ns{\xtn-\xt} 
	\bigg] \notag\\
	& \leq \E\bigg[(1-\theta_t)f(w_{t})+\theta_t f(x^*) 
	-\frac{\mu(1+\mu\eta_t-\alpha_t)\theta_t}{2(1+\mu\eta_t)}\ns{\uxt-x^*} \notag \\
	&\qquad \qquad
	+ \frac{\alpha_t\theta_t}{(1+\mu\eta_t)\eta_t}\Big(\frac{1}{2}\ns{x_{t}-x^*} -\frac{1+\mu\eta_t}{2}\ns{x_{t+1}-x^*}\Big)
	\bigg],  \label{eq:use-eta}
	\end{align}
	where the last inequality \eqref{eq:use-eta} holds by letting $\etat \leq  \frac{\alpha_t}{L(1+\mu\eta_t)\theta_t(1+1/(1-\theta_t))}$.
	
	Finally, according to the probabilistic update of $w_{t+1}$ in Line \ref{line:prob} of Algorithm~\ref{alg:anita}, we have
	\begin{align}
	\E[f(\wtn)] = \E\big[p_t f(\bxtn)+(1-p_t)f(\wt)\big]  \label{eq:use-p}
	\end{align}
	The proof is finished by combining \eqref{eq:use-eta} with \eqref{eq:use-p}, i.e.,  \eqref{eq:keyrelation} is obtained by adding $p_t$ $\times$ \eqref{eq:use-eta} and \eqref{eq:use-p}.
\end{proofof}

\subsection{Proofs for general convex case}
\label{sec:app-1}

In Appendix \ref{sec:app-thm1}, we provide the proof for the main convergence Theorem~\ref{thm:1} in the general convex case (i.e., $\mu=0$). 
Note that the strong convexity Assumption~\ref{asp:strong} is not needed in this case.
Then we provide the proof for its Corollary~\ref{cor:1} with detailed convergence result in Appendix \ref{sec:app-cor1}. 

\subsubsection{Proof of Theorem~\ref{thm:1}}
\label{sec:app-thm1}

First, according to the probabilistic update of $w_{t+1}$ in Line \ref{line:prob} of Algorithm~\ref{alg:anita}, i.e.,
\begin{align}
\wtn= \begin{cases}
\bxtn &\text {with probability } p_t\\
\wt &\text {with probability } 1-p_t	
\end{cases}
\end{align}
Let $p_t\equiv p$ for $0\leq t\leq t_1$, where $t_1$ denotes the first time such that $w_{t_1+1}=\bar x_{t_1+1}$, i.e., $w_{t_1}=w_{t_1-1}=\cdots=w_0$.
Then $t_1\sim \text{Geom}(p)$ and $\E[t_1]=\frac{1-p}{p}$ according to  \eqref{eq:factexp}.
Now, we provide a key Lemma~\ref{lem:first} for the first stage, which shows the decrease of function value in iterations $0\leq t\leq t_1$, and then provide its proof.

\begin{lemma}\label{lem:first}
	Suppose Assumption~\ref{asp:smooth} holds.  
	For $0\leq t\leq t_1$, let $p_t\equiv p$, $\theta_t \equiv \theta$, $\etat \leq  \frac{1}{L(1+1/(1-\theta_t))}$ and $\alpha_t =\theta_t$. 
	Then the following equation holds for \anita (Algorithm~\ref{alg:anita}):
	\begin{align}
		\E[f(w_{t_1+1})-f(x^*)]  
		& \leq 
		\E\Big[(1-\theta) \big(f(x_0)- f(x^*)\big) 
		+ \Big(\frac{\theta^2p}{2\eta}+(1-p)L(1-\theta)\theta^2\Big)\ns{x_0-x^*} \notag\\
		&\qquad \qquad
		-\Big(\frac{\theta^2p}{2\eta}-(1-p)L(1-\theta)\theta^2\Big)\ns{x_{t_1+1}-x^*}\Big].
		\label{eq:first}
	\end{align}
\end{lemma}

\begin{proofof}{Lemma~\ref{lem:first}}
	First, in view of $L$-smoothness of $f$ (Assumption~\ref{asp:smooth}), we recall \eqref{eq:reuse} (where $\forall \beta_t>0$):
	\begin{align}
	&\E[f(\bxtn)] \notag\\
	&\leq \E\bigg[f(\uxt) + \frac{\beta_t}{2L}\ns{\nabla f(\uxt) - \tnabla}  
	+\frac{L(1+1/\beta_t)\theta_t^2}{2}\ns{\xtn-\xt}  
	+ \inner{\nabla f(\uxt)}{\theta_t(x^*-\xt)} \notag \\
	&\qquad \qquad \qquad
	+\inner{\tnabla}{\theta_t(\xtn-x^*)} \bigg] \notag \\
	&= \E\bigg[f(\uxt) 
	+ \frac{\beta_t}{2L}\ns{\nabla f(\uxt) - \tnabla}  
	+\frac{L(1+1/\beta_t)\theta_t^2}{2}\ns{\xtn-\xt} \notag \\
	&\qquad \qquad \qquad
	+ \innerB{\nabla f(\uxt)}{\theta_t(x^*-\uxt)+(1-\theta_t)(\wt-\uxt)}	
	+\inner{\tnabla}{\theta_t(\xtn-x^*)} \bigg] \label{eq:use-uxt1}\\
	&\leq
	\E\bigg[(1-\theta_t)f(\wt) + \theta_t f(x^*) 
	+ \frac{\beta_t}{2L}\ns{\nabla f(\uxt) - \tnabla} 
	+\frac{L(1+1/\beta_t)\theta_t^2}{2}\ns{\xtn-\xt} \notag \\
	&\qquad \qquad \qquad
	+\inner{\tnabla}{\theta_t(\xtn-x^*)}\bigg]    \label{eq:use-convex}\\
	&\overset{\eqref{eq:var1}}{\leq} 
	\E\bigg[(1-\theta_t)f(\wt) + \theta_t f(x^*) 
	+ \frac{L \beta_t }{2}\ns{\uxt-\wt}
	+\frac{L(1+1/\beta_t)\theta_t^2}{2}\ns{\xtn-\xt} \notag \\
	&\qquad \qquad \qquad
	+\inner{\tnabla}{\theta_t(\xtn-x^*)}\bigg] \notag\\
	&=
	\E\bigg[(1-\theta_t)f(\wt) + \theta_t f(x^*) 
	+ \frac{L \beta_t \theta_t^2}{2}\ns{\xt-\wt}
	+\frac{L(1+1/\beta_t)\theta_t^2}{2}\ns{\xtn-\xt} \notag \\
	&\qquad \qquad \qquad
	+\inner{\tnabla}{\theta_t(\xtn-x^*)}\bigg]  \label{eq:use-uxt2}\\
	&=
	\E\bigg[(1-\theta_t)f(\wt) + \theta_t f(x^*) 
	+ \frac{L (1-\theta_t) \theta_t^2}{2}\ns{\xt-\wt}
	+\frac{L(1+1/(1-\theta_t))\theta_t^2}{2}\ns{\xtn-\xt} 
	\notag \\
	&\qquad \qquad \qquad
	+\inner{\tnabla}{\theta_t(\xtn-x^*)}\bigg], \label{eq:first-tmp1}
	\end{align}
	where \eqref{eq:use-uxt1} and \eqref{eq:use-uxt2} use the interpolation step $\uxt=\theta_t\xt+(1-\theta_t)\wt$ (see Line \ref{line:uxt} of Algorithm~\ref{alg:anita}),
	\eqref{eq:use-convex} uses the convexity of $f$, 
	and \eqref{eq:first-tmp1} holds by choosing $\beta_t=1-\theta_t$.
	
	For the last inner product term in \eqref{eq:first-tmp1}, we recall \eqref{eq:tmp-inner} here:
	\begin{align}
	&\E\big[\inner{\tnabla}{\theta_t(\xtn-x^*)}\big] \notag\\
	&\leq \E\bigg[ \frac{\alpha_t\theta_t}{(1+\mu\eta_t)\eta_t} \Big(\frac{1}{2}\ns{\xt-x^*}-\frac{1+\mu\eta_t}{2}\ns{\xtn-x^*} 
	+\frac{\mu\eta_t}{2} \ns{\uxt-x^*}
	-\frac{1}{2}\ns{\xt-\xtn}\Big) \bigg] \notag\\
	&= \E\bigg[ \frac{\alpha_t\theta_t}{\eta_t} \Big(\frac{1}{2}\ns{\xt-x^*}-\frac{1}{2}\ns{\xtn-x^*} 
	-\frac{1}{2}\ns{\xt-\xtn}\Big) \bigg], \label{eq:mu0}
	\end{align}
	where the last equality \eqref{eq:mu0} holds due to $\mu=0$ in this general non-strongly convex case.
	
	Now, we plug \eqref{eq:mu0} into \eqref{eq:first-tmp1} to get 
	\begin{align}
	\E[f(\bxtn)] 
	&\leq
	\E\bigg[(1-\theta_t)f(\wt) + \theta_t f(x^*) 
	+ \frac{\alpha_t\theta_t}{\eta_t} \Big(\frac{1}{2}\ns{\xt-x^*}-\frac{1}{2}\ns{\xtn-x^*} \Big) \notag \\
	&\qquad \qquad
	+ \frac{L (1-\theta_t) \theta_t^2}{2}\ns{\xt-\wt}
	-\frac{\alpha_t\theta_t-L(1+1/(1-\theta_t))\theta_t^2\eta_t}{2\eta_t}\ns{\xt-\xtn}
	\bigg]. \label{eq:first-tmp2}
	\end{align}
	
	According to the parameter setting in Lemma~\ref{lem:first}, we know that $p_t\equiv p$, $\theta_t \equiv \theta$, $\etat \leq  \frac{1}{L(1+1/(1-\theta_t))} \equiv \frac{1}{L(1+1/(1-\theta))} $ and $\alpha_t =\theta_t\equiv \theta$ for iterations $0\leq t\leq t_1$ in the first stage.
	By plugging these parameters into \eqref{eq:first-tmp2}, we obtain 
	\begin{align}
	&\E[f(w_{t_1+1})] \notag\\
	&=\E[f(\bar{x}_{t_1+1})] \notag\\
	&\leq
	\E\bigg[(1-\theta)f(w_0) + \theta f(x^*) 
	+ \frac{\theta^2}{2\eta} \Big(\ns{x_{t_1}-x^*}-\ns{x_{t_1+1}-x^*} \Big) 
	+ \frac{L (1-\theta) \theta^2}{2}\ns{x_{t_1}-w_0}
	\bigg] \notag\\
	&\overset{\eqref{eq:factn1}}{=}
	\E\bigg[(1-\theta)f(w_0) + \theta f(x^*) 
	+ \frac{\theta^2p}{2\eta(1-p)} \Big(\ns{x_{0}-x^*}-\ns{x_{t_1}-x^*} \Big) 
	+ \frac{L (1-\theta) \theta^2}{2}\ns{x_{t_1}-w_0}
	\bigg] \notag \\
	&\overset{\eqref{eq:factn}}{=}
	\E\bigg[(1-\theta)f(w_0) + \theta f(x^*) 
	+ \frac{\theta^2p}{2\eta(1-p)} \Big(\ns{x_{0}-x^*}
	-(1-p)\ns{x_{t_1+1}-x^*} -p\ns{x_{0}-x^*} \Big) 	\notag \\
	&\qquad \qquad
	+ (1-p)\frac{L(1-\theta) \theta^2}{2}\ns{x_{t_1+1}-w_0}
	+p\frac{L(1-\theta) \theta^2}{2}\ns{x_{0}-w_0}
	\bigg] \notag \\
	&\leq 
	\E\bigg[(1-\theta)f(x_0) + \theta f(x^*) 
	+ \frac{\theta^2p}{2\eta(1-p)} \Big(\ns{x_{0}-x^*}
	-(1-p)\ns{x_{t_1+1}-x^*} -p\ns{x_{0}-x^*} \Big) 	\notag \\
	&\qquad \qquad
	+ (1-p)L(1-\theta) \theta^2 \Big(\ns{x_{t_1+1}-x^*}
	+\ns{x_0-x^*}\Big)
	\bigg], \label{eq:use-cauchy} 
	\end{align}
	where \eqref{eq:use-cauchy} uses Cauchy-Schwarz inequality and $w_0=x_0$. 
	Now, the proof of Lemma~\ref{lem:first} is finished since \eqref{eq:first} directly follows from \eqref{eq:use-cauchy}. 
\end{proofof}

After the first stage, for iterations $t>t_1$, we will use a dynamic change of $p_t$. We provide the following Lemma~\ref{lem:t-t1} which shows the decrease of function value in iterations $t> t_1$, and then provide its proof.

\begin{lemma}\label{lem:t-t1}
	Suppose Assumption~\ref{asp:smooth} holds.  
	For $t>t_1$, let $p_t=\max\{\frac{4}{t-t_1+3\sqrt{n}}, \frac{4}{n+3}\}$, $\theta_t =\frac{2}{p_t(t-t_1+3\sqrt{n})}$, $\etat \leq  \frac{1}{3L}$ and $\alpha_t =\theta_t$.
	Then the following equation holds for \anita (Algorithm~\ref{alg:anita}) for any iteration $t>t_1+1$:  
	\begin{align}
		\E\bigg[\frac{\eta_{t-1}}{p_{t-1}\theta_{t-1}^2}\big(f(w_{t})-f(x^*)\big)\bigg]  
		& \leq \E\bigg[\frac{(1-p_{t_1+1}\theta_{t_1+1})\eta_{t_1+1}}{p_{t_1+1}\theta_{t_1+1}^2}\big(f(w_{t_1+1})-f(x^*)\big)
		\notag\\
		&\qquad \qquad  \qquad
		+ \frac{1}{2}\Big(\ns{x_{t_1+1}-x^*} -\ns{x_{t}-x^*} \Big) \bigg].
		\label{eq:t-t1}
	\end{align}
\end{lemma}

\begin{proofof}{Lemma~\ref{lem:t-t1}}
	For proving this lemma, we will use our technical Lemma~\ref{lem:func_relation} with $\mu=0$ (general convex case).
	In particular, by choosing $\alpha_t=\theta_t$ and multiplying $\frac{\eta_t}{p_t\theta_t^2}$ for both sides in \eqref{eq:keyrelation} with $\mu=0$, we obtain
    the following Lemma~\ref{lem:change}:
	
	\begin{lemma}\label{lem:change}
		Suppose Assumption~\ref{asp:smooth} holds.  
		Choose stepsize $\etat \leq  \frac{1}{L(1+1/(1-\theta_t))}$ and $\alpha_t =\theta_t$ for any $t\geq0$.
		Then the following equation holds for \anita (Algorithm~\ref{alg:anita}) for any iteration $t\geq 0$:
		\begin{align}
			\E\bigg[\frac{\eta_t}{p_t\theta_t^2}\big(f(w_{t+1})-f(x^*)\big)\bigg]  
			\leq \E\bigg[\frac{(1-p_t\theta_t)\eta_t}{p_t\theta_t^2}\big(f(w_{t})-f(x^*)\big)
			+ \frac{1}{2}\Big(\ns{x_t-x^*} -\ns{x_{t+1}-x^*} \Big) \bigg].
			\label{eq:change}
		\end{align}
	\end{lemma}

	Then, we are going to sum up \eqref{eq:change} from iteration $t_1+1$ to $t$ for obtaining \eqref{eq:t-t1}.
	In order to get a recursion formula for \eqref{eq:change}, we further choose appropriate parameters $\{p_t\}$, $\{\theta_t\}$ and $\{\eta_t\}$ to obtain 
	\begin{align}
	\frac{(1-p_t\theta_t)\eta_t}{p_t\theta_t^2} \leq \frac{\eta_{t-1}}{p_{t-1}\theta_{t-1}^2}. \label{eq:pt-thetat}
	\end{align}
	It is not hard to verify that \eqref{eq:pt-thetat} can be satisfied for any $t>t_1+1$ by choosing 
	\begin{align}
	&p_t=\max\Big\{\frac{4}{t-t_1+3\sqrt{n}}, \frac{4}{n+3}\Big\}, \label{eq:pt} \\
	&\theta_t =\frac{2}{p_t(t-t_1+3\sqrt{n})}, \label{eq:thetat} \\
	&\eta_{t} \equiv \eta \leq \frac{1}{3L},
	\end{align}
	for any $t> t_1$.
	The proof of Lemma~\ref{lem:t-t1} is finished by summing up \eqref{eq:change} from iteration $t_1+1$ to $t$ and noting that \eqref{eq:pt-thetat} holds for any $t>t_1+1$.
\end{proofof}

Now, we are ready to prove Theorem~\ref{thm:1} by combining Lemma~\ref{lem:first} (iterations $0\leq t\leq t_1$) and Lemma~\ref{lem:t-t1} (iterations $t>t_1$).

\begin{proofof}{Theorem~\ref{thm:1}}
	First, we note that $p_{t_1+1}=\frac{4}{1+3\sqrt{n}}$, $\theta_{t_1+1}=\frac{1}{2}$ and $\eta_{t_1+1} \leq  \frac{1}{L(1+1/(1-\theta_{t_1+1}))}=\frac{1}{3L}$.
	Then, by plugging \eqref{eq:first} into \eqref{eq:t-t1} and noting that
	\begin{align}
	\frac{(1-p_{t_1+1}\theta_{t_1+1})\eta_{t_1+1}}{p_{t_1+1}\theta_{t_1+1}^2} 
	= \frac{(1-\frac{2}{1+3\sqrt{n}})\frac{1}{3L}}{\frac{1}{1+3\sqrt{n}}}
	=\frac{3\sqrt{n}-1}{3L} 
	\leq \frac{4\sqrt{n}}{L},  \label{eq:mt1}
	\end{align}
	we have
	\begin{align}
	&\E\bigg[\frac{\eta_{t-1}}{p_{t-1}\theta_{t-1}^2}\big(f(w_{t})-f(x^*)\big)\bigg] \notag\\  
	& \leq \E\bigg[\frac{4\sqrt{n}}{L} \bigg((1-\theta) \big(f(x_0)- f(x^*)\big) 
	+ \Big(\frac{\theta^2p}{2\eta}+(1-p)L(1-\theta)\theta^2\Big)\ns{x_0-x^*}  \notag\\
	&\qquad \qquad   
	-\Big(\frac{\theta^2p}{2\eta}-(1-p)L(1-\theta)\theta^2\Big)\ns{x_{t_1+1}-x^*}\bigg) 
	+ \frac{1}{2}\Big(\ns{x_{t_1+1}-x^*} -\ns{x_{t}-x^*} \Big) \bigg] \notag\\
	&\leq \E\bigg[\frac{4\sqrt{n}}{L} \bigg((1-\theta) \big(f(x_0)- f(x^*)\big) 
	+ \Big(\frac{\theta^2p}{2\eta}+(1-p)L(1-\theta)\theta^2\Big)\ns{x_0-x^*}\bigg) \notag\\
	&\qquad \qquad   
	-\bigg(\Big(\frac{\theta^2p}{2\eta}-(1-p)L(1-\theta)\theta^2\Big)\frac{4\sqrt{n}}{L}  -\frac{1}{2}\bigg) \ns{x_{t_1+1}-x^*}  \bigg] \notag\\
	&\leq \E\bigg[\frac{4\sqrt{n}}{L} \bigg((1-\theta) \frac{L}{2}
	+\frac{\theta^2p}{2\eta}+(1-p)L(1-\theta)\theta^2\bigg)\ns{x_0-x^*} \notag\\
	&\qquad \qquad   
	-\bigg(\Big(\frac{\theta^2p}{2\eta}-(1-p)L(1-\theta)\theta^2\Big)\frac{4\sqrt{n}}{L}  -\frac{1}{2}\bigg) \ns{x_{t_1+1}-x^*}  \bigg] \label{eq:f0}\\
	&\leq 8\ns{x_0-x^*}, \label{eq:removet1}
	\end{align}
	where \eqref{eq:f0} holds due to the $L$-smoothness of $f$,
	\eqref{eq:removet1} follows from the constant parameters  i.e., $p_t\equiv p=\frac{1}{n+1}$, $\theta_t\equiv \theta=1-\frac{1}{2\sqrt{n}}$ and  $\eta_{t} \leq  \frac{1}{L(1+1/(1-\theta_{t}))}\equiv \frac{1}{(1+2\sqrt{n})L}$ for $t\leq t_1$.
	
	Finally, the proof of Theorem~\ref{thm:1} is finished by multiplying $\frac{p_{t-1}\theta_{t-1}^2}{\eta_{t-1}}$ for both sides of \eqref{eq:removet1}, i.e.,  we have for any iteration $t>t_1+1$
	\begin{align}
	\E[f(w_{t})-f(x^*)]  
	&\leq \frac{8p_{t-1}\theta_{t-1}^2\ns{x_0-x^*} }{\eta_{t-1}} \overset{\eqref{eq:thetat}}{=}
	\frac{32\ns{x_0-x^*}}{\eta_{t-1}p_{t-1}(t-t_1+3\sqrt{n})^2}. \label{eq:final1}
	\end{align}
\end{proofof}

\subsubsection{Proof of Corollary~\ref{cor:1}}
\label{sec:app-cor1}
Now, we provide the proof for Corollary~\ref{cor:1} with detailed convergence result of \anita in the general convex case (i.e., $\mu=0$). 

\begin{proofof}{Corollary~\ref{cor:1}}
	Note that the output of \anita (Algorithm~\ref{alg:anita}) is $w_T$ after $T$ iterations.
	To show that $w_T$ is an $\epsilon$-approximate solution, we recall \eqref{eq:final1} with iteration $t=T>t_1+1$ here:
	\begin{align}
	\E[f(w_{T})-f(x^*)]  \leq
	\frac{32\ns{x_0-x^*}}{\eta_{T-1}p_{T-1}(T-t_1+3\sqrt{n})^2}. \label{eq:cor-tmp1}
	\end{align}
	According to \eqref{eq:pt}, we know $p_t=\max\Big\{\frac{4}{t-t_1+3\sqrt{n}}, \frac{4}{n+3}\Big\}$ for any $t>t_1$. 
	Thus we divide \eqref{eq:cor-tmp1} into two cases, i) $p_t=\frac{4}{t-t_1+3\sqrt{n}}$ for $t_1<t\leq t_1+n+3-3\sqrt{n}$; 
	ii) $p_t=\frac{4}{n+3}$ for $t>t_1+n+3-3\sqrt{n}$.
	
	Now, we know that for Case i) $t\leq t_1+n+3-3\sqrt{n}$,  then $p_t=\frac{4}{t-t_1+3\sqrt{n}}$, $\theta_t=\frac{2}{p_t(t-t_1+3\sqrt{n})}=\frac{1}{2}$, $\eta_{t} \leq \frac{1}{3L}$, and \eqref{eq:cor-tmp1} turns to 
	\begin{align}
	\E[f(w_{T})-f(x^*)]  \leq
	\frac{24L\ns{x_0-x^*}}{T-t_1+3\sqrt{n}} \leq \epsilon. \label{eq:cor-tmp2}
	\end{align}
	The last inequality of \eqref{eq:cor-tmp2} holds by choosing $T=t_1-3\sqrt{n} +\frac{24L\ns{x_0-x^*}}{\epsilon}$.
	In particular, if $\epsilon\geq O(\frac{1}{n})$, then (recall that $\E[t_1]=n$ and also it can be derandomized to $n$ iterations)
	\begin{align}
	T=t_1-3\sqrt{n} +\frac{96L\ns{x_0-x^*}}{\epsilon}\leq 2n. \label{eq:tstage1}
	\end{align}
	
	For the other case $\epsilon< O(\frac{1}{n})$ (small convergence error), it corresponds to Case ii) $t> t_1+n+3-3\sqrt{n}$ (i.e., more iterations are needed),  then 
	$p_t=\frac{4}{n+3}$, $\theta_t=\frac{2}{p_t(t-t_1+3\sqrt{n})}=\frac{n+3}{2(t-t_1+3\sqrt{n})}\leq \frac{1}{2}$,  $\eta_{t} \leq \frac{1}{3L}$ and \eqref{eq:cor-tmp1} turns to 
	\begin{align}
	\E[f(w_{T})-f(x^*)]  \leq
	\frac{24(n+3)L\ns{x_0-x^*}}{(T-t_1+3\sqrt{n})^2} \leq \epsilon. \label{eq:cor-tmp3}
	\end{align}
	The last inequality of \eqref{eq:cor-tmp3} holds by choosing
	\begin{align}
	T=t_1-3\sqrt{n}+\sqrt{\frac{24(n+3)L\ns{x_0-x^*}}{\epsilon}}
	\leq n+\sqrt{\frac{24(n+3)L\ns{x_0-x^*}}{\epsilon}}. \label{eq:tstage2}
	\end{align}
	
	Now, the remaining thing is to bound the number of stochastic gradient computations of \anita for achieving the $\epsilon$-approximate solution $w_T$.
	As we discussed in Section \ref{sec:alg}, we know that \anita (Algorithm~\ref{alg:anita})  uses $(n+2)p_{t}+2(1-p_{t})=np_t+2$ stochastic gradients in expectation for iteration $t$. 
	According to the choice of probability $\{p_t\}$ in Corollary~\ref{cor:1} (Theorem~\ref{thm:1}), we know that there are three stages.
	1) The first stage $p_t\equiv\frac{1}{n+1}$ for $0\leq t\leq t_1$;  2) the second stage $p_t=\frac{4}{t-t_1+3\sqrt{n}}$ for $t_1<t \leq t_1+n +3-3\sqrt{n}$; 3) the third stage $p_t\equiv \frac{4}{n+3}$ for $t>t_1+n +3-3\sqrt{n}$.
	
	First, let us consider the case of large $\epsilon$ (i.e., $\epsilon\geq O(\frac{1}{n})$).
	Then we know that only the first two stages of \anita is enough for finding an $\epsilon$-approximate solution in this case.
	According to \eqref{eq:cor-tmp2}, the total number of stochastic gradient computations is
	\begin{align}
	\#\mathrm{grad} 
	=\sum_{t=0}^{T-1}(np_t+2)
	&=n\Big(\sum_{t=0}^{t_1}p_t + \sum_{t=t_1+1}^{T-1}p_t\Big)+ 2T \notag\\
	&= n\Big(\sum_{t=0}^{t_1}\frac{1}{n+1} + \sum_{t=t_1+1}^{T-1}\frac{4}{t-t_1+3\sqrt{n}}\Big) + 2T \notag\\
	&\leq O\bigg(n\Big(1+\log\frac{1}{\epsilon\sqrt{n}}\Big)\bigg), \label{eq:gradcase1}
	\end{align}
	where the last inequality \eqref{eq:gradcase1} follows from \eqref{eq:tstage1}.
	
	Then, for the other case $\epsilon< O(\frac{1}{n})$ (small convergence error), we know that more iterations are needed for finding an $\epsilon$-approximate solution.
	According to \eqref{eq:cor-tmp3}, the total number of stochastic gradient computations is
	\begin{align}
	\#\mathrm{grad} 
	&=\sum_{t=0}^{T-1}(np_t+2) \notag\\
	&=n\Big(\sum_{t=0}^{t_1}\frac{1}{n+1} ~~
	+ ~ \sum_{t=t_1+1}^{t_1+n+3-3\sqrt{n}}\frac{4}{t-t_1+3\sqrt{n}} ~~
	+ ~ \sum_{t=t_1+n+4-3\sqrt{n}}^{T-1}\frac{4}{n+3}\Big)+ 2T \notag\\
	&\leq O\left(n\log\sqrt{n}+\sqrt{\frac{nL}{\epsilon}}\right), \label{eq:gradcase2}
	\end{align}
	where the last inequality \eqref{eq:gradcase2} follows from \eqref{eq:tstage2}.
\end{proofof}

\subsection{Proofs for strongly convex case}
\label{sec:app-strong}

Similar to Appendix \ref{sec:app-1}, we first provide the proof of the main convergence Theorem~\ref{thm:strong} for the strongly convex case (i.e., $\mu>0$) in Appendix \ref{sec:app-thmstrong}. 
Then we provide the proof for its Corollary~\ref{cor:strong} with detailed convergence result in Appendix \ref{sec:app-corstrong}. 

\subsubsection{Proof of Theorem~\ref{thm:strong}}
\label{sec:app-thmstrong}

As we discussed at the end Remark of Section \ref{sec:strong},
the parameter setting of \anita in this strongly convex case is simpler than the general convex case in Theorem~\ref{thm:1}. 
As a result, we only need one technical Lemma~\ref{lem:strong} in this proof sketch of Theorem~\ref{thm:strong} rather than three Lemmas \ref{lem:first}--\ref{lem:t-t1} in previous general convex case.
Note that Lemma~\ref{lem:strong} directly follows from our technical Lemma~\ref{lem:func_relation} with $\alpha_t=1+\mu\eta_{t}$.

\begin{lemma}\label{lem:strong}
	Suppose that Assumptions \ref{asp:smooth} and \ref{asp:strong} hold. 
	Choose stepsize $\etat \leq  \frac{1}{L\theta_t(1+1/(1-\theta_t))}$ and $\alpha_t =1+\mu\eta_t$ for any $t\geq0$. 
	Then the following equation holds for \anita (Algorithm~\ref{alg:anita}) for any iteration $t\geq 0$:  	
	\begin{align}
		&\E\bigg[f(w_{t+1})-f(x^*) + \frac{(1+\mu\eta_t)p_t\theta_t}{2\eta_t}\ns{x_{t+1}-x^*}\bigg] 
		\leq \E\bigg[(1-p_t\theta_t)\big(f(w_{t})-f(x^*)\big) + \frac{p_t\theta_t}{2\eta_t}\ns{x_{t}-x^*}\bigg]. 
		\label{eq:strong-1}
	\end{align}
\end{lemma}

\begin{proofof}{Theorem~\ref{thm:strong}}
	According to Lemma~\ref{lem:strong}, we know that the change of function value after a gradient update step in \anita. Then, according to the parameter settings chosen in Theorem~\ref{thm:strong}, we have $p_t\equiv p$ and $\theta_t \equiv \theta=\frac{1}{2}\min\{1,\sqrt{\frac{\mu}{pL}}\}$ for any $t\geq0$, and the  stepsize $\etat \leq  \frac{1}{L\theta_t(1+1/(1-\theta_t))} \equiv \eta= \frac{1}{L\theta(1+1/(1-\theta))}$.
	Now, we further define 
	\begin{align}
	\Phi_t:=f(w_t)-f(x^*)+\frac{(1+\mu\eta)p\theta}{2\eta}\ns{x_t-x^*},  \label{eq:defphi}
	\end{align}
	then \eqref{eq:strong-1} in Lemma~\ref{lem:strong} can be changed to, for any iteration $t\geq 0$, 
	\begin{align}
	\E[\Phi_{t+1}]  
	& \leq \E\bigg[\max\Big\{1-p\theta,\ \frac{1}{1+\mu\eta}\Big\} \Phi_t\bigg] \notag\\
	&\leq \E\bigg[\Big(1-\frac{4p\theta}{5}\Big) \Phi_t\bigg] \label{eq:use-theta-eta} \\
	&\leq \Big(1-\frac{4p\theta}{5}\Big)^{t+1} \Phi_0. \label{eq:final-strong}
	\end{align}
	where \eqref{eq:use-theta-eta} uses $\frac{1}{1+\mu\eta} \leq 1-\frac{4p\theta}{5}$ since the choice of parameters $\theta=\frac{1}{2}\min\{1,\sqrt{\frac{\mu}{pL}}\}$ and $\eta= \frac{1}{L\theta(1+1/(1-\theta))}$, and the last inequality \eqref{eq:final-strong} holds by telescoping \eqref{eq:use-theta-eta} from iteration $t$ to $0$.
\end{proofof}

\subsubsection{Proof of Corollary~\ref{cor:strong}}
\label{sec:app-corstrong}
Now, we provide the proof for Corollary~\ref{cor:strong} with detailed convergence result of \anita in the strongly convex case (i.e., $\mu>0$). 

\begin{proofof}{Corollary~\ref{cor:strong}}
	Note that the output of \anita (Algorithm~\ref{alg:anita}) is $w_T$ after $T$ iterations.
	To show that $w_T$ is an $\epsilon$-approximate solution, we recall \eqref{eq:final-strong} with iteration $t=T-1$:
	\begin{align}
	\E[f(w_{T})-f(x^*)]  
	\leq \E[\Phi_{T}] 
	\leq
	\Big(1-\frac{4p\theta}{5}\Big)^T \Phi_0\leq \epsilon, \label{eq:cor-tmp-strong}
	\end{align}
	where the first inequality is due to the definition of $\Phi_T$ (see \eqref{eq:defphi}), and the last inequality holds by letting the number of iterations 
	$$T= \frac{5}{4p\theta}\log\frac{\Phi_0}{\epsilon}.$$ 
	Moreover, by choosing $p=\frac{1}{n}$ and recalling that $\theta=\frac{1}{2}\min\{1,\sqrt{\frac{\mu}{pL}}\}$, then the total number of stochastic gradient computations of \anita for achieving the $\epsilon$-approximate solution $w_T$ is
	\begin{align*}
	\#\mathrm{grad} 
	=\sum_{t=0}^{T-1}(np_t+2)
	=\Big(n\frac{1}{n}+2\Big)T
	=\frac{15}{4p\theta}\log\frac{\Phi_0}{\epsilon}
	=O\left(\max\Bigg\{n, \sqrt{\frac{nL}{\mu}}\Bigg\}\log\frac{1}{\epsilon}\right).
	\end{align*}
\end{proofof}

\end{document}